\documentclass[a4paper]{article}

\usepackage{graphicx}
\usepackage{psfrag}
\usepackage[english]{babel}
\usepackage{geometry}
\usepackage{amssymb}
\usepackage{amsfonts}
\usepackage{amsmath}
\usepackage{amscd}
\usepackage{amsthm}
\usepackage[all]{xy}
\usepackage{color}
\usepackage{bbm}
\usepackage{cite}

\newtheoremstyle{plainNoItalics}{}{}{\normalfont}{}{\bfseries}{.}{ }{}

\newtheorem{theo}{Theorem}[section]
\newtheorem{prop}[theo]{Proposition}

\newtheorem{lem}[theo]{Lemma}

\newtheorem*{theo*}{Theorem \ref{theorem-main}}

\newtheorem*{theorem*}{Theorem}
\newtheorem*{proposition*}{Proposition}
\newtheorem*{lemma*}{Lemma}
\newtheorem*{corollary*}{Corollary}

\theoremstyle{plainNoItalics}
\newtheorem{defi}[theo]{Definition}
\newtheorem*{definition*}{Definition}
\newtheorem{rem}[theo]{Remark}
\newtheorem{ex}[theo]{Example}

\newtheorem{conjecture}[theo]{Conjecture}

\newtheorem*{remark*}{Remark}

\newtheorem*{observation*}{Observation}
\newtheorem*{example*}{Example}
\newtheorem*{assumption*}{Assumption}

\input xy
\xyoption{all}

\begin{document}
\begin{center}
{\huge Conley Index at Infinity} \\
\vspace{1cm}
Juliette Hell\\
Freie Universit\"at Berlin, Germany\\
blanca@math.fu-berlin.de\\
\vspace{0.5cm}
\begin{minipage}{9cm}
\small{ {\bf Abstract:} The aim of this paper is to explore the possibilities of Conley index techniques in the study of heteroclinic connections between finite and infinite invariant sets. For this, we remind the reader of the Poincar\'e compactification: this transformation allows to project a $n$-dimensional  vector space $X$ on the $n$-dimensional unit hemisphere of $X\times \mathbb{R}$ and infinity on its $(n-1)$-dimensional equator called the sphere at infinity. Under a  normalizability condition, vector fields on $X$ are mapped to vector fields on the  Poincar\'e hemisphere whose associated flows leave the equator invariant. The dynamics on the equator reflects the dynamics at infinity, but are now finite and may be studied by Conley index techniques. Furthermore, we observe that some non-isolated behavior may occur around the equator,  and introduce the concept of an  invariant set at infinity  of isolated invariant dynamical  complement. Through the construction of an extended phase space together with an extended flow, we are able to adapt the Conley index techniques and prove the existence of connections to such non-isolated invariant  sets.}
\end{minipage}
\end{center}
\psfrag{Xx1}{$X\times\{ 1 \}$}
\psfrag{H}{$\mathcal{H}$}
\psfrag{E}{$\mathcal{E}$}
\psfrag{PMm}{$\mathcal{P}(x)$}
\psfrag{Mm}{$(x,1)$}
\psfrag{O}{$O$}
\psfrag{z}{$z$}
\psfrag{00}{$r_2$}
\psfrag{20}{$s_1$}
\psfrag{-20}{$s_2$}
\psfrag{r1}{$r_1$}
\psfrag{a1}{$a_1$}
\psfrag{a2}{$a_2$}
\psfrag{18}{18}
\psfrag{25}{25}
\psfrag{39}{39}
\psfrag{A}{}
\psfrag{B}{}
\psfrag{Scomp}{$S_{comp}$}
\psfrag{S}{$S$}
\psfrag{B1}{$B_1$}
\psfrag{B2}{$B_2$}
\psfrag{r}{$r$}
\psfrag{a}{$a$}
\psfrag{c}{$c$}
\psfrag{B+}{$B^+$}
\psfrag{B-}{$B^-$}
\psfrag{b+}{$b^+$}
\psfrag{b-}{$b^-$}

\section{Introduction}

The analysis of the structure of connections in global attractors inspired several tools, in particular the Conley index techniques. Phenomena like blow up or grow up may be interpreted as heteroclinic connections to infinity. Examples for such situations are to be found in \cite{messias, kappos2}\\
In this paper, we describe a way of using the Conley index theory for the analysis of unbounded  global attractors - or more generally unbounded invariant sets. For this, we need to analyse both the finite connections and those going to infinity in forward or in backward time. \\
The first step of our method consists in making infinity finite by Poincar\'e compactification. The phase space is projected to a hemisphere of the same dimension, infinity being projected onto its equator.  The equator is the boundary of the hemisphere and is itself a sphere  called the sphere at infinity. Its dimension is (dimension of the phase space)$-1$. For more details see Section \ref{comp} and in particular Figure \ref{projpoinc}. This procedure grants us two things:  the exhibition of the relevant dynamic at infinity, i.e. the invariant sets in the sphere at infinity, and the consideration of {\bf bounded} neighbourhoods of these.\\
However, Conley index does not only require bounded neighborhoods of invariant sets, but isolating ones. For more details on classical Conley index theory, see the original work \cite{conley} and good introductions in \cite{smollerbook, kappos3, kappos4}. In the question of isolation, it turns out  that infinity, although compactified, still stands out by virtue of  its lack of isolation in many cases. To illustrate this, let us have a look at Figure \ref{classhamilton}: it shows  equilibria at infinity that are not isolated, in phase portraits 18 and 39. \\
To circumvent this difficulty, we introduce the notion of a complementary isolated invariant set. The precise definition will be given in \ref{defcomp}. Let us here illustrate the three possible situations with three examples.  We pick up those examples in \cite{schlovu, arllib} where classifications of (quadratic) planar vector fields are given.  The numbers associated to those phase portrait coincide with their numbers in the papers. 
\begin{enumerate}
\item  The phase Portrait 25 of Figure (\ref{classhamilton}) is from \cite{arllib}. It admits only isolated invariant equilibria at infinity.
\item The phase Portrait 18 of Figure (\ref{classhamilton}) is from \cite{arllib}. It admits an equilibrium  at the bottom which is not isolated invariant but of isolated invariant dynamical complement.
\item The phase Portrait 39 of  Figure (\ref{classhamilton}) comes from \cite{schlovu}. It admits an equilibrium  at the bottom  which is neither isolated invariant nor complementarily isolated.  
\end{enumerate}
A Conley index at infinity for such invariant sets is defined and denoted by $\hat{h}(.)$. This index is based on time-duality and the classical index of the dynamical complement - an interesting new concept that we introduce in Definition \ref{defcomp} .With its help, one can state basic existence/nonexistence theorems of the following type (see Theorem (\ref{existence}) for more details, as well as Definition \ref{indexofattrep}). 
\begin{theorem*}
Consider an invariant set $S$ at infinity  of isolated invariant complement. If the  Conley index at infinity of $S$, $\hat{h}(S)$,  is neither of  a repeller nor of an attractor, then there are trajectories accumulating on $S$ in forward time direction, and trajectories accumulating on $S$ in backward time direction. 
\end{theorem*}
Furthermore we introduce a construction in order to analyse the connection structure of an unbounded invariant set, provided it consists of isolated invariant parts or of invariant parts of isolated invariant complement. The rough idea is to replace the parts that are invariant and only complementarily isolated by an ersatz infinity. This will be done in Section \ref{constrersatz}.
After proceeding with this construction, we are able to prove the main Theorem \ref{detect}. The content of this theorem can be summarized as follows.
\begin{theorem*}
If the  Conley index theory proves the existence of connecting orbits from/to the isolated invariant  ersatz infinity to/from an isolated invariant set $R\subset S_{comp}$ for the extension $\hat{\mathcal{H}}$ with extended flow $\hat{\phi}$, \\
then the corresponding connecting orbits $S\rightarrow R$/$R\rightarrow S$ exist for the original compactified flow $\phi$ on  the Poincar\'e Hemisphere $\mathcal{H}$.
\end{theorem*}
Hence we are able to detect heteroclinic orbits to sets which, a priori, were beyond the reach of Conley index theory. Such invariant sets pop up naturally when studying transfinite dynamics and the connection structure of unbounded invariant sets. 

The author thanks the referee for his constructive remarks, and  Ahmad Afuni for his careful corrections reading. This research has been supported by the Deutsche Forschungsgemeinschaft, SFB 647-Space-Time-Matter.

\section{Poincar\'e compactification} \label{comp}
The Poincar\'e compactification transforms a finite dimensional vector space $X$ into a compact manifold -- a hemisphere -- and is therefore called a compactification. The construction also works  for infinite dimensional Hilbert spaces $X$, in which case the image would be  an infinite dimensional hemisphere, which is a noncompact Hilbert manifold with boundary. Arbitrarily far points in the vector space are mapped onto the equator -- the boundary of the Poincar\'e hemisphere. Therefore we will also refer to it as the sphere at infinity. The sphere at infinity is a sphere whose dimension reads $\dim (X)-1$.
We refer to the overview paper \cite{kappos2}, and to the books \cite{perko, ALGM} where explicit planar examples are computed. Now, let us describe  the  construction of the Poincar\'e compactification more precisely.

\begin{figure}
\begin{center}
\includegraphics[width=0.8\textwidth]{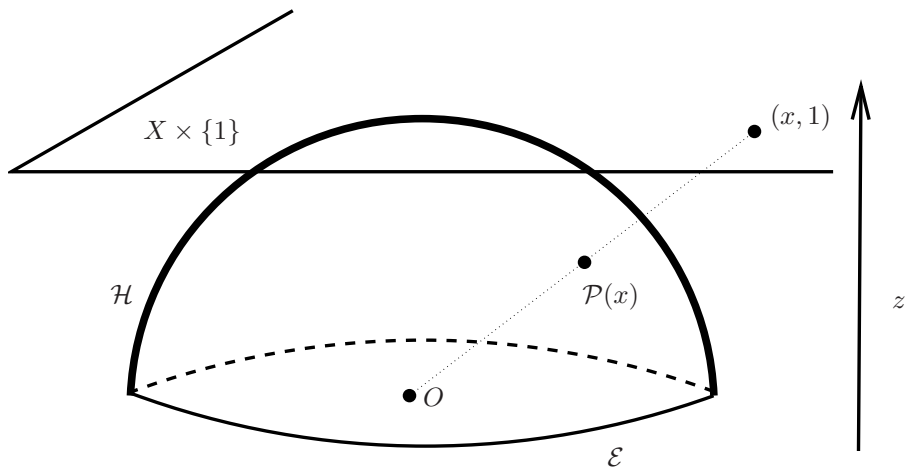}
\end{center}
\caption{ \label{projpoinc}  The Poincar\'e transformation.} 
\end{figure}

Let $X$ be an Hilbert space whose scalar product is denoted by $\left\langle  \cdot ,\cdot  \right\rangle $. We consider the vector space $X\times \mathbb{R}$. In other words, we add to $X$ a, say, vertical dimension whose coordinate will be denoted by $z\in \mathbb{R}$. The vector space $X\times  \mathbb{R}$ is naturally  equipped with the scalar product $\left\langle (x,z), (x',z') \right\rangle := \left\langle   x,x' \right\rangle +zz'$. The original vector space $X$ is identified with the affine hyperplane $X \times \{1\}$ through the bijection $x\in X \mapsto (x,1)\in X\times \{1\}$. Finally, the affine hyperplane $X\times \{1\}$ is projected centrally to the unit hemisphere 
\[ \mathcal{H}:=\{(x,z) \in X \times \mathbb{R} / \langle  x,x \rangle  + z^2 =1, z\geqslant 0 \},  \]
called the Poincar\'e hemisphere. 
This transformation is sketched in Figure (\ref{projpoinc}).

The Poincar\'e transformation is  given by the formula
\begin{eqnarray*} 
 \mathcal{P}:& X \rightarrow &\mathcal{H} \\
             & x \mapsto     & (\chi,z) = \frac{1}{   \sqrt{ \left\langle x,x\right\rangle +1} } (x,1).
\end{eqnarray*}
As $\Vert x \Vert$ tends to infinity, $\mathcal{P}(x)$ goes to the equator 
\[\mathcal{E}:=\{ (\chi,0) \in X\times \mathbb{R} / \left\langle \chi , \chi  \right\rangle =1  \}\]
of the Poincar\'e hemisphere. Hence the equator $\mathcal{E}$ is also called the sphere at infinity.

Now let us describe how vector fields on the vector space $X$ are mapped under the Poincar\'e transformation. For this purpose, let us consider the equation 
\begin{equation} \label{basiceqn} 
\dot{x}=f(x), \quad x\in X,
\end{equation}
where the upper  dot denotes time derivative.
Applying the Poincar\'e transformation, we get  the following evolution equations for the image $(\chi, z)=\mathcal{P}(x)$, $z>0$:
 \begin{equation}\label{basicchizeqn}
 \left\{
 \begin{array}{ccc}
\dot{\chi} &=& f_{z} (\chi) - \left\langle f_{z} (\chi), \chi \right\rangle \chi, \\
\dot{z}&=&- \left\langle  f_{z} (\chi), \chi \right\rangle z
\end{array}
\right.
\end{equation}
where $f_z$  denotes the homothetic of $f$ with factor $z$; in other words
\[f_z(\chi):= z f( z^{-1} \chi).\]
Equation (\ref{basicchizeqn}) makes sense on the sphere at infinity if and only if $f_z$ has a limit as $z$ tends to zero. This is rarely the case. 
However this difficulty can be overcome by normalization: we  replace  $f_z$ by $\rho(z)f_z$ in Equation (\ref{basicchizeqn}), where $\rho$ is 
\begin{enumerate}
\item a positive function,
\item strictly positive on $\mathcal{H}\setminus{\mathcal{E}}$ so that the finite trajectories are not affected, 
\item tends quickly enough to zero as $z\rightarrow 0$, slowing down the trajectories and guaranteeing the existence of a continuous limit: 
\end{enumerate}
\[\lim_{z \rightarrow 0} \rho( z )f_{z}(\chi)=f_0(\chi).\]
More precisely, normalizable vector fields are defined as follows.
\begin{defi}{\bf Normalization of the homothety $f_z$}\label{normpoincc}\\
The vector field $f:X\rightarrow X$ is normalizable if and only if there exists a continuous function $\rho: [0,+\infty [ \rightarrow [0,+ \infty [$ such that \begin{enumerate}
\item $\rho(z)\neq 0$ if $z\neq 0$, and 
\item the map $(\chi, z)\mapsto \rho(z)f_z(\chi)$ is Lipschitz continuous even at $z=0$.
\end{enumerate} 
\end{defi}
\begin{rem} {\bf Normalization is a change of time variable.} \\
Replacing the homothety $f_z({\chi})$ by its normalized version $\rho(z)f_z({\chi})$ corresponds to a change of time variable $dt=\rho(z)d\tau$.
We will nevertheless keep the name $t$ for the time variable of the normalized equations on the Poincar\'e hemisphere for the sake of simplicity. When we speak of a flow/vector field on the Poincar\'e hemisphere, we mean such a normalized equation where regularity all the way up to the equator is given. 
\end{rem}
\begin{rem}
{\bf Vector fields with polynomial growth are normalizable.}
 Consider a map $f:X \rightarrow X$ admitting a decomposition of the following type: 
\[f=P+p,\]
where $P$ is homogenous of degree $d$ ( i.~e.~ $P(\lambda x)= \lambda^d P(x)$)
and $p$ contains terms of order lower than $d-1$, but not necessarily polynomial  (i.~e.~ $p(x)    =   o(\Vert x \Vert ^d)$ for $\Vert x \Vert$ going to infinity). 
The normalization by $\rho(z)=z^{d-1}$ provides 
\[\lim_{z \rightarrow 0}\rho(z) f_z = P.\]
\end{rem}
We will from now on mostly omit  mention of  normalization, and speak only of compactified vector field/flow on the Poincar\'e hemisphere $\mathcal{H}$, meaning in fact compactified and normalized. Otherwise we do not even get a proper flow on $\mathcal{H}$. 

At last, we give some more formulas which are useful for computations. Instead of analyzing a  vector field near the equator of the Poincar\'e hemisphere (typically, looking for its equilibria  at infinity and their stability), it is easier to compute in vertical hyperplanes tangent to the hemisphere  around the region of interest. Therefore, we need the following projection $Proj_e$ onto the vertical affine half hyperplane $E$ tangent to $\mathcal{H}$ at $e \in \mathcal{E}$, i.~e.~ $E=\left( e+<e>^{\perp} \right) \cap \left( X \times [0,\infty[   \right)$. We call  the coordinates in these planes the ``vertical charts". 
\begin{eqnarray*}
Proj_e  :& \mathcal{H}& \rightarrow E \\
               &  (\chi,z )   & \mapsto \frac{1}{\left\langle  \chi, e \right\rangle} (\chi, z) :=(\xi, \zeta)
\end{eqnarray*}
The resulting vector field on $E$ in the variables $(\xi, \zeta)$ reads
\begin{equation} \label{evolxi}
\left\lbrace 
\begin{array}{lll}
\dot{\xi}& = & -\left< f_{\zeta}(\xi), e \right>  \xi +  f_{\zeta}(\xi)  \\
\dot{\zeta} & = & - \left<  f_{\zeta}(\xi), e \right> \zeta 
\end{array}
\right.
\end{equation}
Here again, it is often necessary to replace $f_{\zeta}$ by its normalized form $\rho f_{\zeta}$ to realize regularity of the projected vector field on the vertical charts.

After compactification, the sphere at infinity may contain isolated invariant sets whose Conley indices are well defined. We recall the definition of the Conley index for an isolated invariant set intersecting the boundary of the Poincar\'e hemisphere: it splits into three components because the standard definition is not sufficient to take the dynamics on the  boundary $\partial \mathcal{H}= \mathcal{E}$ into consideration. This definition coincides with those of \cite{mccord, mrozek}.\\
Let us first fix  some vocabulary. The property of being open - or closed - is to be understood relative to the Poincar\'e hemisphere $\mathcal{H}$. We write, when possible, a subscript $\mathcal{H}$ to remind the reader of this. Furthermore, a set $K$ is said to be a compact neighborhood if $K$ is compact and the closure of its interior $K=cl_{\mathcal{H}}(int_{\mathcal{H}}(K))$. An (forward and backward) invariant set $S$ is said to be isolated invariant iff there is a compact neighborhood $K$ of $S$ for which $S$ is the maximal invariant set $Inv(K)=S$ and $S\subset int_{\mathcal{H}}(K)$. 
\begin{defi} \label{conleyboundary}  {\bf Conley index on the boundary.}\\
We consider an isolated invariant set $S$ on the Poincar\'e hemisphere  $\mathcal{H}$ with boundary $\mathcal{E}$. The set $S$ admits an index pair $(N,N_1)$, i.e. a pair of compact set satisfying 
\begin{enumerate}
\item {\bf Isolation:} The set  $cl(N\setminus{N_1})$ is an isolating neighbourhood for  $S$. More precisely, $S\subset cl(N\setminus{N_1})$ does not intersect the boundary of $cl_{\mathcal{H}}(N\setminus{N_1})$ relatively to $\mathcal{H}$.  
\item {\bf Positive invariance of $N_1$ with respect to $N$:} A trajectory starting in $N_1$ remains in $N_1$ until it leaves $N$, i.~e. for all $x\in N_1$ and all $t>0$, $\varphi ([0,t], x)\subseteq N \Longrightarrow \varphi ([0,t], x)\subseteq N_1$.
\item {\bf The set $N_1$ is an exit set for $N$: } The trajectories leave $N$  through $N_1$, i.~e. for all $x\in N$ and $t_{1}>0$ with $\varphi(t_{1}, x )\notin N$, there exists a time $t_{0}\in [0,t_{1}]$,  such that $\varphi ([0,t_{0}], x )\subseteq N$ and $\varphi (t_{0} ,x)\in N_1$ .
\end{enumerate}
Let $*\notin X$ be a universal point. The Conley index of $S$ is defined by three homotopy classes of pointed spaces as follows:\\
$h(\mathcal{H};S):=\left[ \frac{N\cup \{*\}}{N_1\cup \{*\} } \right]$, the Conley index with respect to $\mathcal{H}$;\\
 $h(\mathcal{H}, \mathcal{E};S):=\left[ \frac{N\cup \{*\}}{N_1 \cup (N \cap \mathcal{E})\cup \{*\}} \right]$, the Conley index with respect to $\mathcal{H}$ and $\mathcal{E}$;\\
$h(\mathcal{E};S):=\left[ \frac{N \cap \mathcal{E} \cup \{*\}}{N_1\cap \mathcal{E} \cup \{*\}} \right]$, the Conley index with respect to $\mathcal{E}$.
\end{defi}

The standard Conley index is able to analyse the dynamics at infinity in some situations, as illustrated in the following. 
\begin{ex} 
\label{perkoindices}
{\bf A quadratic vector field under Poincar\'e compactification}\\
We consider the following polynomial system in the plane: 
\begin{equation} \label{257}
\left\lbrace 
\begin{array}{lll}
\dot{x_1} &=& x_1^2+x_2^2-1\\
\dot{x_2} &=& 5(x_1x_2-1)
\end{array}
\right.
\end{equation}
This system does not show any finite equilibria, hence also no finite periodic orbit. The dynamic  at infinity was analysed in the vertical charts, and Figure (\ref{perko257_ex_poinc}) gives an overview of the global phase portrait on the Poincar\'e hemisphere, seen from above. All equilibria at infinity are isolated invariant and their indices can be easily computed and read as follows:
\begin{equation*}
 \begin{array}{ll}
h(a_1) =h(a_2)=\left\lbrace   
                              \begin{array}{lll} 
                                                            h(\mathcal{H};a_{1,2 }) &= &\Sigma^0 \\
                                                            h(\mathcal{H}, \mathcal{E} ;a_{1,2 } ) &= & \bar{0}\\
                                                            h( \mathcal{E} ; a_{1,2 }) &=& \Sigma^0 
                               \end{array}
         \right.
&
  h(r_1)= h(r_2)= \left\lbrace   
                              \begin{array}{lll} 
                              h(\mathcal{H}; r_{1,2}) &= & \bar{0}\\
                              h(\mathcal{H}, \mathcal{E} ; r_{1,2}) &= & \Sigma^2 \\
                              h( \mathcal{E} ; r_{1,2}) &= & \Sigma^1
                              \end{array}
           \right.
 \\
 h(s_1) = \left\lbrace   
                              \begin{array}{lll} 
                              h(\mathcal{H};s_1 ) &= & \bar{0}\\
                              h(\mathcal{H}, \mathcal{E} ;s_1 ) &= & \Sigma^1 \\
                              h( \mathcal{E} ; s_1) &= &  \Sigma^0
                              \end{array}
           \right.
 &
  h(s_2) = \left\lbrace   
                              \begin{array}{lll} 
                              h(\mathcal{H};s_2 ) &= & \Sigma^1\\
                              h(\mathcal{H}, \mathcal{E} ; s_2) &= &\bar{0} \\
                              h( \mathcal{E} ; s_2) &= & \Sigma^1
                              \end{array}
           \right.
 \end{array}
\end{equation*}
%
%
%

\begin{figure}[h]
\begin{center}
\includegraphics[width=0.4\textwidth]{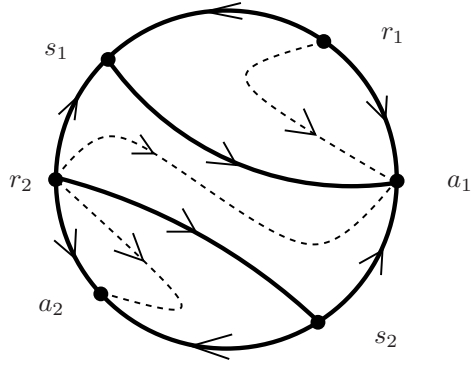}
\caption{The Poincar\'e compactification of the vector Field (\ref{257}).\label{perko257_ex_poinc}} 
\end{center}

\end{figure}
\end{ex}

\begin{figure}
\begin{center}
\includegraphics[width=1\textwidth]{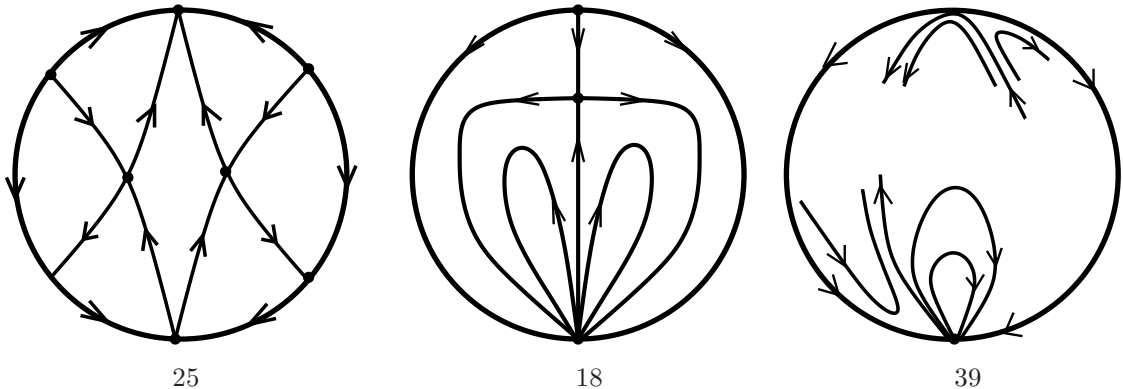}
\end{center}
\caption{ \label{classhamilton}  Equilibria at infinity.} 
\end{figure}

On the other hand, isolated invariance seems to be too strong a requirement for the dynamics at infinity.  Let us briefly motivate the next section through  examples. Quadratic vector fields in the plane may exhibit 
accumulations of homoclinic loops at infinity as illustrated in the phase portrait 18, 39 of  Figure (\ref{classhamilton}). 
Quadratic planar vector fields have been studied and classified:\\
In \cite{arllib} quadratic planar Hamiltonian vector fields are classified according to their full phase portraits on the Poincar\'e hemisphere. There are 28 different phase portraits, 21 of them exhibit only isolated invariant equilibria at infinity, while 7 of them contain equilibria at infinity which are not isolated invariant because of an accumulation of homoclinic loops. Those 7 still fit our  definition of  complementary isolated invariance introduced in Definition \ref{defcomp}. \\
In \cite{schlovu}, quadratic planar vector fields are classified according to their phase portraits near the sphere at infinity. This classification contains 40 phase portraits: 25 of them show only isolated invariant equilibria at infinity, 12 of them show  equilibria of isolated invariant complement, and  3 of them   show equilibria which are neither isolated invariant nor of isolated invariant complement. \\
Figure \ref{classhamilton} illustrates these three cases. The Phase Portraits 25, 18 are taken from \cite{arllib}: Portrait  25 exhibits only isolated invariant equilibria; Portrait 18 admits an equilibrium of isolated invariant complement at the bottom. The Phase Portrait 39 comes from \cite{schlovu} and admits  an equilibrium at the bottom which is neither isolated invariant nor of isolated invariant dynamical complement.    

For a polynomial vector field in finite dimensions, it seems that an accumulation of homoclinic loops on an equilibrium at infinity  is structurally stable, as soon as this equilibrium is hyperbolic with respect to the flow within the sphere at infinity (as for instance in Figure (\ref{classhamilton}), Portrait  18). For more details, see Theorem \ref{strucstab} in Paragraph \ref{genericity} which was proven in \cite{gonzvel}. To be able to deal with some of these situations, we develop in the following paragraph the concept of an invariant set of isolated invariant dynamical complement and the Conley index at infinity. 
\section{Conley index at infinity}
\subsection{Isolation of the complement}
The existence of an accumulation of homoclinic loops at infinity even for planar quadratic fields motivates the following definition of an invariant set of isolated invariant dynamical complement. The basic idea is the following. If it is not  possible to isolate an invariant set $S$, try to isolate the invariant set $S_{comp}$ containing every trajectory that does not tend to the first set $S$, neither in backward nor in forward time direction. In other words, $S_{comp}$ contains the trajectories that remain far away from $S$ and can therefore be considered as its dynamical complement. The precise definition is given in Definition (\ref{defcomp}). This idea is comparable to the concept of an  attractor repeller pair, but here we do not  require the connecting  trajectories to run on a one-way-street.  The compact neighborhoods isolating  such an isolated invariant set $S_{comp}$   are to be chosen arbitrarily large, as long as they do not intersect $S$ itself.

The complement of a set $K$ is denoted by $K^c=\mathcal{H}\setminus{K}$, not to be confused with the dynamical complement of a set $S$ denoted by $S_{comp}$ and defined in Definition (\ref{defcomp}).  

\begin{defi} \label{defcomp}  {\bf dynamical complement}.\\
Consider a normalized  flow on the Poincar\'{e} Hemisphere $\mathcal{H}$.
Let $S$ be a  closed (forward and backward) invariant set. We call  the set $S_{comp}$ 
\[  S_{comp}:= \{ x\in \mathcal{H} :\quad  \alpha(x)\cap S = \omega(x) \cap S = \emptyset   \}\]
the  dynamical complement of $S$.
Furthermore the invariant set $S$ is said to be of isolated invariant dynamical complement if $S_{comp}$ is isolated invariant. 
\end{defi}
The set $S_{comp}$ is invariant because for every $y\in o(x)$, $\alpha(x)=\alpha(y)$ and $\omega(x)= \omega(y)$. 

\begin{ex} \label{18}
Let us have a closer look at the phase portrait 18 of Figure (\ref{classhamilton}). The equilibrium $S$ at infinity at the bottom is not isolated because of the accumulation of homoclinic loops.  The equilibrium $S$ is  of isolated invariant complement. Its dynamical complement $S_{comp}$ consists in 
\begin{enumerate}
\item The repeller at infinity at the top, 
\item the saddle at the origin,
\item the trajectory connecting the first to the second. 
\end{enumerate}
The dynamical complement $S_{comp}$ is invariant and isolated by every compact neighbourhood of itself  that does not contain  the equilibrium $S$.\\
\end{ex}

Determining the dynamical complement $S_{comp}$ from $S$ requires much knowledge on the global dynamics, so that one may want to formulate equivalent definitions of  complementary isolated invariance where $S_{comp}$ itself does not come into play.
The first equivalent definition requires isolating properties for every large  compact neighborhood that does not intersect $S$.  The second one only for a continuous family of them, which is more easily verifiable  in concrete examples. 
\begin{prop} \label{defequiv1}
We consider a normalized flow on the Poincar\'e hemisphere. A closed invariant set $S$ is of isolated invariant complement if and only if  there exists a compact neighbourhood $K\subset \mathcal{H}$ with the following properties.
\begin{enumerate}
\item The set $K$ does not intersect  the invariant set $S$: 
\[   K\cap S = \emptyset .\]
\item  The set $K$ is an isolating neighbourhood in the Poincar\'e hemisphere:
\[  Inv(K)\subset int_{\mathcal{H}} (K).  \]
\item Every compact neighbourhood $K' \supset K$ which does not intersect the invariant set  $S$ is also an isolating neighbourhood:
\begin{equation*}
\left.
\begin{array}{r}
K'  \text{ compact neighbourhood}\\
K' \supset K \\
K' \cap S = \emptyset
\end{array}
\right\rbrace
\Rightarrow
K' \text{ isolating neighbourhood.}
\end{equation*}
\end{enumerate} 
\end{prop}
\begin{proof}
Suppose that the dynamical complement $S_{comp}$ of $S$ is isolated invariant. Hence $S_{comp}$ is closed as a maximal invariant set,   disjoint from $S$,  and there exists an isolating neighborhood $K$ of $S_{comp}$ that does not intersect $S$. We claim that this $K$ also satisfies  the third condition of Proposition 
\ref{defequiv1}. Indeed, let $K'$ be a compact neighborhood containing $K$ with $K'\cap S=\emptyset$, and $x\in Inv(K')$. Then the orbit through $x$ as well as its $\alpha-$ and $\omega-$limit sets, is contained in $K'$ and therefore do not intersect $S$. By definition, such an $x$ is in $S_{comp}$. Hence $Inv(K')\subset S_{comp}$. As $K'\supset K$,  the reverse inclusion $Inv(K')\supset Inv(K)=S_{comp}$ holds , and the claim $Inv(K')=S_{comp}$ is proved. \\
Now let us assume the existence of a compact neighborhood $K$ satisfying the three conditions of Proposition \ref{defequiv1}.  First of all we prove that $K$ as well as every superset $K'$ as in condition (3) isolate the same invariant set. We proceed indirectly. If there exists a $K'$ satisfying condition (3) with $Inv(K')\varsupsetneq Inv(K)$, then there is a point $x\in Inv(K')\setminus{ Inv(K)}$, whose orbit $o(x)$ is fully contained in $K'$ but leaves $K$. 
Consider $\tilde{K}=K \cup k$, where $k$ is a compact neighborhood with $\left( o(x)\setminus{K} \right)\cap \partial k \neq \emptyset$ and $k\subset K'$ - for example half of  a small enough tubular neighborhood around $o(x)\setminus{K}$. Then  $o(x)\subset Inv(\tilde{K})$ holds, but $o(x)\cap \tilde{K}\neq \emptyset$, so that $\tilde{K}\subset K$ cannot be an isolating neighborhood which is a contradiction to condition (3).\\
Now let us prove that $Inv(K)=S_{comp}$. As $K\cap S=\emptyset$, the inclusion $Inv(K)\subset S_{comp}$  surely holds true. Now consider any point $x\in S_{comp}$. If the point $x$ were not in $Inv(K)$, then it  would also not belong to any $Inv(K')$ for $K'$ of condition (3), because we just proved $Inv(K)=Inv(K')$. Taking a sequence $K_n$ of such sets so that $d(K_n, S)\rightarrow 0$, we get a sequence of points $x_n\in o(x)$ with the property
\[ 
x_n \notin K_n, \quad d(x_n,S)\rightarrow 0.
\] 
As $\mathcal{H}$ is compact, the sequence $(x_n)$ converges (up to subsequence) to a point $\xi \in S$. The point $\xi$ belongs to the orbit of $x$, or  to its $\alpha$- or $\omega$-limit set, contradicting $x\in S_{comp}$. Hence $Inv(K)=S_{comp}$. 
\end{proof}
\begin{ex}
Let us illustrate this proposition with the phase portraits of Figure (\ref{classhamilton}). 

Portrait 39: The equilibrium at the top is not of isolated invariant complement. If $K^c$ is a small enough neighborhood of the upper equilibrium, then $K$ is not an isolating neighborhood: there are trajectories arbitrarily near this equilibrium that produce internal tangencies to $\partial_{\mathcal{H}}K$ but never leave $K$, hence $Inv(K)\not\subset int_{\mathcal{H}}(K)$. 

Portrait 18: The bottom equilibrium is of isolated invariant complement. If $K$ is any {\it big enough} compact neighborhood that does not contain the bottom equilibrium, then every point of $\partial_{\mathcal{H}}K$ leaves $K$ in forward or backward time direction. Hence $Inv(K)\subset int_{\mathcal{H}}(K)$ i.e. the compact set $K$ is an isolating neighborhood.  Choosing $K$ to small leads to violation of condition (3) as seen in Remark \ref{SandVhom}. 

\end{ex}

The third equivalence is the most easily verifiable in concrete examples. It ``only" requires a continuous family of isolating neighbourhoods whose complements shrink on the complementary isolated invariant set $S$. But in fact it turns out to be a {\it stronger} requirement than that in Proposition (\ref{defequiv1}), as we will explain in Remark (\ref{SandVhom}), so that we have to exclude some pathologic cases to keep the equivalence.
\begin{defi} \label{wild}  {\bf wild.} \\
We call a set $S$ wild if it does not admit any compact neighborhood $V$ such that $S$ is a deformation retract of $V$. 
\end{defi}
\begin{rem}
This definition of wild may not be completely standard, but this is the condition we need to guarantee the direction $\Rightarrow$ in the equivalence below. An hawaiian earring is wild in the sense of Definition \ref{wild}.
\end{rem}
\begin{prop} \label{def2}
We consider a  flow $\phi$ on the Poincar\'e hemisphere.
Let be $S$ an invariant set that is not wild. Then
$S$ is  of isolated invariant complement iff there exists a compact neighborhood  $V$ of $S$ and a deformation retraction 
\[  H: [0,1]\times V \rightarrow V \]
from $V$ to $S$ such that for every $\lambda \in [0,1[$, the compact set $K_{\lambda}:= cl_{\mathcal{H}}\left( H(\lambda, V)^c \right)$ is an isolating neighborhood in the Poincar\'e hemisphere $\mathcal{H}$. 
The map $H$ being a deformation retraction means:
\begin{eqnarray*}
\forall x \in V, &  H(0,x)=x ,& \\
\forall x \in V, & H(1,x) \in S,& \\
\forall x \in S, & H(1, x)=x ,& \\
\end{eqnarray*}
\end{prop}
\begin{rem}\label{SandVhom}
The fact that $S$ is a deformation retract of its neighborhood $V$ implies in particular that they are both of the same homotopy type. It is not possible for wild sets, therefore  we have to exclude them. To be more precise, let us consider an invariant hawaiian earring in phase portrait 18: $S$ shall now be the union of the bottom equilibrium with an infinite but countable number of homoclinic loops of size tending to 0. This set $S$ is an invariant set of isolated invariant complement $S_{comp}$ (the same $S_{comp}$ as in Example \ref{18}) for which we cannot find any continuous family of isolating neighborhoods as defined in Proposition \ref{def2} because the invariant set $S$ considered is wild. \\
If we do not assume that $S$ is not wild in Proposition \ref{def2}, then we only have 
\[ \text{Existence of a continuous family }K_{\lambda} \text{ as in \ref{def2}} \Rightarrow S_{comp}  \text{ isolated invariant}\]
On the other hand, the requirement is strong enough to prevent the modification of the invariant set isolated by the $K_{\lambda}$ all along the deformation. Of course we will prove this below, but let us consider again $S=$bottom equilibrium of phase portrait (18), and refer again to Example \ref{18}. Suppose we consider a continuous family of isolating neighborhoods with $K_0=$union of two disjoint  squares around the repeller at the top and the saddle at the origin (hence yet to small to isolate the whole $S_{comp}$). As the $K_{\lambda}$'s grow, the two squares collide and form a rectangle around the orbit connecting the top and origin equilibria. The invariant set isolated by the continuous family changes  without violating the isolation property. This does not contradict Proposition \ref{def2}, because this continuous family is {\it not} admissible there: $S$ is not a deformation retract of $V:=cl_{\mathcal{H}}(K_0^c)$, because $V$ is not simply connected.  

\end{rem}
\begin{proof}
Assume the existence of a continuous family of isolating neighborhoods. 
First we prove that for all $\lambda \in [0,1)$, the sets $K_{\lambda}= cl_{\mathcal{H}}(H(\lambda, V)^c)$ isolate the same maximal invariant set. For that proceed indirectly: if there were an orbit $o(x)$ through $x\in K_0$ and a $\tilde{\lambda}>0$ with
\begin{eqnarray*}
 && o(x) \subset K_{\tilde{\lambda}}, \\
o(x) \nsubseteq K_0 &\Leftrightarrow& o(x)\cap V \neq \emptyset
\end{eqnarray*}
The set $\{ \lambda \in [0,\tilde{\lambda}] / \quad o(x)\cap H(\lambda, V)\neq \emptyset \}$ contains 0 and admit $\tilde{\lambda }$ as an upper bound, so that it admits a supremum
\begin{equation}\label{lambda1sup} \lambda_1 := \sup \{ \lambda \in [0,\tilde{\lambda}]  / \quad o(x)\cap H(\lambda, V)\neq \emptyset \}, \qquad 0\leq \lambda_1 < \tilde{\lambda}
\end{equation}
Claim 1: $\begin{cases} o(x)\subset K_{\lambda_1} \\ o(x) \cap \partial K_{\lambda_1} \neq \emptyset \end{cases}$, where $K_{\lambda_1}=cl_{\mathcal{H}}(H(\lambda_1, V)^c)$, hence $K_{\lambda_1}$ cannot be an isolating neighborhood.\\
If $o(x)\nsubseteq K_{\lambda_1}$, then there exists a point $y\in int_{\mathcal{H}}(H(\lambda_1, V))\cap o(x)$. But then $y\in int_{\mathcal{H}}(H(\lambda_2, V))\cap o(x)$ for $\lambda_2$ slightly bigger than $\lambda_1$, contradiciting the maximality of $\lambda_1$ (\ref{lambda1sup}). Hence $o(x)\subset K_{\lambda_1}$.\\
 Now let us show $ o(x) \cap \partial K_{\lambda_1} \neq \emptyset$. By definition (\ref{lambda1sup}) of $\lambda_1$, there exists a sequence $(\mu_k)_{k \in \mathbb{N}}\to \lambda_1$ and sequences of points $(x_k)_{k \in \mathbb{N}}$, $(y_k)_{k \in \mathbb{N}}$ with
 \[  x_k = H(\mu_k, y_k) \in o(x)\cap H(\mu_k, V), \quad y_k \in V. \]
 The sequence  $(y_k)_{k \in \mathbb{N}}$ converges by compacity of $V$, up to a subsequence, to a point $y\in V$. By continuity of $H$, we have
 \[ H(\mu_k, y_k)=x_k \in o(x) \longrightarrow  H(\lambda_1, y) \in H(\lambda_1, V), \qquad \text{ as }k\to \infty . \]
 This provides a point $H(\lambda_1, y)\in o(x) \cap \left( K_{\lambda_1} \cap H(\lambda_1, V) \right)= o(x) \cap \partial K_{\lambda_1}$, as claimed.\\
Now we prove that the common isolated invariant set $T$ to all $K_{\lambda}$, $\lambda \in [0,1)$, is the dynamical complement $S_{comp}$ of $S$. First, $T\subset S_{comp}$ because $T\subset V^c$ implies $\alpha(x)\cap S = \omega(x) \cap S = \emptyset$ for every $x\in T$. Furthermore proceed again indirectly:  if $S_{comp}\nsubseteq T = Inv(K_{\lambda})$ for all $\lambda \in [0,1)$, then 
\begin{eqnarray*}
\forall \lambda \in [0,1),& \quad \exists& x_\lambda \in S_{comp}\\
&&o(x_\lambda) \cap H(\lambda, V) \neq \emptyset.
\end{eqnarray*} 
Considering a sequence $\lambda \to 1$, we get as a limit a nontrivial intersection $S\cap S_{comp}$, which is a contradiction. Hence $S_{comp}= T = Inv(K_{\lambda})$ for all $\lambda \in [0,1)$, is isolated invariant. \\
Now we prove the second part of the equivalence. Assume that $S_{comp}$ is isolated invariant. By Proposition \ref{defequiv1}, there exists a compact neighborhood $K$ isolating $S_{comp}$ with the property (3) of Proposition \ref{defequiv1}. Because $S$ is not wild, there is a deformation retraction $H$ of a compact neighborhood $V\subset K^c$ to $S$ which provides a continuous family of isolating neighborhoods $K_{\lambda}:= cl_{\mathcal{H}}(H(\lambda, V)^c$, $\lambda \in [0,1)$, isolating $S_{comp}$. 
\end{proof}
\begin{defi}\label{indexatinf} {\bf Conley index at infinity.}\\
Consider a flow $\phi$ on the Poincar\'e hemisphere, and  an invariant set $S\subset \mathcal{E}$  at infinity of isolated invariant complement. Its dynamical complement $S_{comp}$ is well defined and isolated invariant.  We define the Conley index at infinity of $S$, denoted by $\hat{h}(S)$, as the classical Conley index of $S_{comp}$ computed under the flow $\phi_-(t,.):=\phi(-t,.)$ with reversed time direction. 
\end{defi}
\begin{rem}
Let us justify why this provides a well defined index. 

A compact set $R\subset \mathcal{H}$ is isolated invariant under $\phi$ if and only if $R$ is isolated invariant under $\phi_-$. Furthermore, If $B$ is an isolating block for $R$ with respect to the flow $\phi$, and $B_{\pm}$ its immediate entrance and exit sets respectively, then $B$ is also an isolating block for $R$ with respect to the flow $\phi_-$ with respective entrance and exit sets $B_{\mp}$.    By classical Conley index theory, we know
\begin{itemize}
\item  
$(B, B_+) $ is an index pair for $R$ with respect to the flow $\phi_-$, 
\item 
the Conley index $h_{\phi_-}(R)$ of $R$ under the flow $\phi_-$ is defined by the quotients 
\[h_{\phi_-}(R): \begin{cases}
h{\phi_-}(\mathcal{H}; R):=\left[  \frac{B}{B_+}\right] \\ h{\phi_-}(\mathcal{H}, \mathcal{E}; R):=\left[  \frac{B}{B_+ \cup \left(  B\cap \mathcal{E} \right)  }\right] \\  h{\phi_-}( \mathcal{E}; R):=\left[  \frac{B \cap \mathcal{E}}{B_+ \cap \mathcal{E}   }\right] 
\end{cases}
\]
\item and  this definition does not depend on the choice of the index pair. 
\end{itemize}
By Definition \ref{indexatinf} of the Conley index at infinity, we have 
\[\hat{h}_{\phi}(S):=h_{\phi_-}(S_{comp}).\]
The Conley index at infinity of $S$ does not depend on the choice of the index pair for $S_{comp}$.
\end{rem}
\begin{ex}
The indices in Example (\ref{18}), Figure (\ref{classhamilton}), portrait 18,  are the following:
\begin{itemize}
\item The repeller $r$ at infinity at the top  has index 
\[
h(r) =  \begin{cases}
                                  h(\mathcal{H}; r )= \bar{0} \\
                                  h(\mathcal{H}, \mathcal{E}; r )= \Sigma^2 \\
                                   h( \mathcal{E}; r )= \Sigma^1
         \end{cases}.
\]
\item The saddle $O$ at the origin has index $h(O)=\Sigma^1$.
\item The isolated invariant set $S_{comp}$ has Conley index 
\[
h(S_{comp}) =  \begin{cases}
                                  h(\mathcal{H}; S_{comp} )= \Sigma^1 \\
                                  h(\mathcal{H}, \mathcal{E}; S_{comp})= \bar{0} \\
                                   h( \mathcal{E}; S_{comp})= \Sigma^1
         \end{cases}.
\]
\item The Conley index at infinity of the degenerate equilibrium $S$ reads
\[
\hat{h}(S) =  \begin{cases}
                                  \hat{h}(\mathcal{H}; S )= \bar{0} \\
                                  \hat{h}(\mathcal{H}, \mathcal{E}; S )= \Sigma^1 \\
                                   \hat{h}( \mathcal{E}; S )= \Sigma^0
         \end{cases}.
\]
\end{itemize} 
\end{ex}
\begin{rem}
The author is convinced that the concept of invariant set of isolated invariant complement,  together with the concept of  dynamical complement, are  new and crucial for the study of dynamics at infinity or more generally on bounded manifolds with boundary. However, it is not clear wether the Definition (\ref{indexatinf}) is the most sensible for the index of  a set with isolated invariant complement.   An argument against it is the following. If $S$ is both isolated invariant and of isolated invariant complement, then both $h(S)$ and $\hat{h}(S)$ are defined, but we see in Example (\ref{pech})  that they do not have to coincide - basically because the summability formula fails for $\hat{h}$. \\
We want  to point out here that this technical point is not crucial for the detection of heteroclinics to infinity -- which is the goal that we pursue.

Definition (\ref{indexatinf}) is   useful for results on existence/nonexistence of connections to infinity. Such theorems are given in the next paragraph dealing with general properties of invariant sets of isolated invariant complement and their indices. 
\end{rem}
\begin{ex} \label{pech}
We consider a flow on the  3-dimensional Poincar\'e hemisphere, which we represent as a 3-dimensional ball whose boundary is the sphere at infinity. This flow is sketched in Figure (\ref{notcoinc}).
\begin{figure}
\begin{center}
\includegraphics[width=0.4\textwidth]{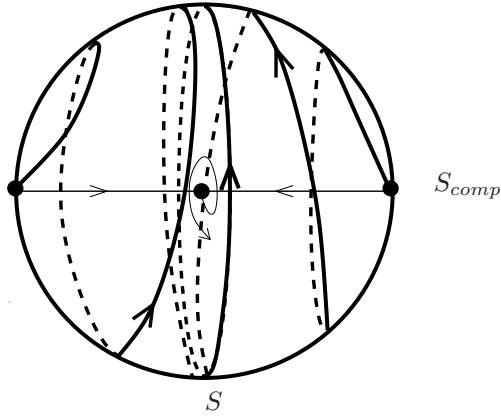}
\end{center}
\caption{ \label{notcoinc} 3-dimensional Poincar\'e hemisphere.} 
\end{figure}
The finite dynamics exhibit only an equilibrium at the origin which is a saddle. The 1-dimensional stable manifold of the origin connects  to two fully unstable equilibria at infinity. The 2-dimensional unstable manifold of the origin contains trajectories swirling toward a periodic orbit at infinity which is fully stable. Furthermore, on the sphere at infinity the two unstable equilibria  are also connected to the stable periodic orbit  by swirling trajectories. \\
Let us call  the periodic orbit in the sphere at infinity $S$. The invariant set $S$ is isolated invariant. Its dynamical complement  $S_{comp}$ is well defined, being the union of the three equilibria and the stable manifold of the origin. Moreover, $S_{comp}$ happens to be isolated invariant too. Therefore we can compare $h(S)$ and $\hat{h}(S)$ in this example. 
A block $B$ for $S_{comp}$ is given by a cylinder around the unstable manifold of the origin connecting the three equilibria. The immediate exit set $B^-$ consists of the whole boundary $B^-=\partial B$. Therefore
$$
h(S_{comp})= \begin{cases}
                          h(\mathcal{H}; S_{comp})= \bar{0} \\
                           h(\mathcal{H, \mathcal{E}}; S_{comp})= \Sigma^3\\
                            h(\mathcal{E}; S_{comp})=  \Sigma^2  \vee  \Sigma^2
                          \end{cases}
\hat{h}(S)=h_-(S_{comp})= \begin{cases}
                          h(\mathcal{H};  S)= \Sigma^0 \\
                           h(\mathcal{H, \mathcal{E}}; S)= \bar{0}\\
                            h(\mathcal{E}; S)=  \Sigma^0  \vee  \Sigma^0.
                          \end{cases}
$$
Now let us consider the classical Conley index $h(S)$ of the periodic orbit at infinity. A block $D$  for $S$ is given by a 3-dimensional ring around $S$, intersected with the Poincar\'e hemisphere. The periodic orbit $S$ being stable, the immediate exit set of the block $D$ is empty and
$$
h(S)= \begin{cases}
                          h(\mathcal{H};  S)= \bigcirc * \\
                           h(\mathcal{H, \mathcal{E}}; S_{comp})= \bar{0}\\
                            h(\mathcal{E}; S)=  \bigcirc *,
                          \end{cases}
$$
where $\bigcirc *$ denotes the disjoint union of a 1-dimensional circle and the universal point $*$, not to be confused with $\Sigma^1$ where the union is not disjoint.  \\
Hence $h(S)\neq \hat{h}(S)$. The topologies of the sets $S$ and $S_{comp}$ are different and this fact is reflected in the indices $h(S)$ and $\hat{h}(S)$.
\end{ex}

\subsection{General properties}
The Conley index at infinity is of course also defined on the homology and cohomology levels as the homology/cohomology of the pointed spaces given by $\hat{h}$. The Poincar\'e duality translates in terms of time duality and provides the following relationships between the indices. The proof of this proposition is a direct application of the result of \cite{mrozek}.
\begin{prop}
Consider a  flow on the $n$-dimensional Poincar\'e hemisphere $\mathcal{H}$. We denote by $H_*$ and $H^*$ the homological and cohomological Conley indices, and by $\hat{H}_*$ and $\hat{H}^*$ the homological and cohomological Conley indices at infinity, respectively. If $S$ is of isolated invariant complement $S_{comp}$, 
\begin{eqnarray*}
\hat{H}^k(\mathcal{H}, \mathcal{E};S) &=& H_{n-k}(\mathcal{H};S_{comp}) \\
\hat{H}^k(\mathcal{H};S) &=& H_{n-k}(\mathcal{H}, \mathcal{E};S_{comp}) \\
\hat{H}^k(\mathcal{E};S) &=& H_{n-k}( \mathcal{E};S_{comp})
\end{eqnarray*}
\end{prop}

Furthermore, the index at infinity allows us to formulate existence or nonexistence results. For this, we describe which types of  Conley indices characterize attractors or repellers (Definition \ref{indexofattrep}) in order to formulate the existence/nonexistence Theorem \ref{existence}.

\begin{defi} \label{indexofattrep}
Consider a flow on the Poincar\'e hemisphere $\mathcal{H}$. We refer to an index  coming from an index pair of the form $(B,\emptyset)$  as an ``index of an attractor" an index . Such an index takes the form:
\begin{eqnarray*}
h(\mathcal{H};S) &=& [B]\cup \{ *\} \\
h(\mathcal{H}, \mathcal{E};S) &=& \bar{0} \\
h(\mathcal{E};S) &=& [B\cap \mathcal{E}]\cup \{ *\} 
\end{eqnarray*} 
We refer to an index coming from an index pair of the form $(B,\partial B)$ as an ``index of a repeller" an index . Such an index takes the form:
\begin{eqnarray*}
h(\mathcal{H};S) &=&  \bar{0} \\
h(\mathcal{H}, \mathcal{E};S) &=& \left[ \frac{B}{\partial_{\mathcal{B}} B}  \right]  \\
h(\mathcal{E};S) &=&\left[ \frac{B\cap \mathcal{E}}{\partial_{\mathcal{H}} B \cap \mathcal{E}}  \right]
\end{eqnarray*} 
\end{defi}
If an invariant set $S$ is complementary isolated and its dynamical complement $S_{comp}$ has the index $h(S_{comp})$ of an attractor or of a repeller, then $S$ itself has a Conley index at infinity $\hat{h}(S)$ of a repeller or of an attractor respectively, due to the time reversal in the definition. 

This concept allows us to formulate the following existence/nonexistence theorem.
\begin{theo}\label{existence}

Consider a complementary isolated invariant set $S\subset \mathcal{E}$ in the sphere at infinity. Its dynamical complement $S_{comp}$ is isolated invariant and $\hat{h}(S)$ is well defined. 
\begin{enumerate}
\item If $\hat{h}(S)$ is the  index of  a repeller, then \begin{itemize}
                                                                                        \item there are trajectories $S\rightarrow S_{comp}$, or more precisely $$\exists x \in \mathcal{H}\setminus{\left( S\cup S_{comp}\right) },\alpha(x)\cap S \neq \emptyset  \text{, and } \omega(x)\subset S_{comp}$$,
                                                                                        \item there is no trajectory $S_{comp}\rightarrow S$.
                                                                                        \end{itemize}
\item  If $\hat{h}(S)$ is the index of  an attractor, then \begin{itemize}
                                                                                        \item there are trajectories $S_{comp}\rightarrow S$, or more precisely $$\exists x \in \mathcal{H}\setminus{\left( S\cup S_{comp}\right) },\alpha(x)\subset S_{comp} \text{, and }\omega(x)\cap S \neq \emptyset  $$,
                                                                                        \item there is no trajectory $S\rightarrow S_{comp}$.
                                                                                        \end{itemize}
\item If $\hat{h}(S)$ is neither the index of an attractor nor of a repeller, then there are both trajectories $S\rightarrow S_{comp}$ and $S_{comp} \rightarrow S$, or more precisely
 $$\exists x \in \mathcal{H}\setminus{\left( S\cup S_{comp}\right) },\alpha(x)\cap S \neq \emptyset  \text{, and } \omega(x)\subset S_{comp}$$  $$\exists y \in \mathcal{H}\setminus{\left( S\cup S_{comp}\right) },\alpha(y)\subset S_{comp} \text{, and }\omega(y)\cap S \neq \emptyset  $$
\end{enumerate}
 
\end{theo}
\begin{proof}
Suppose $\hat{h}(S)$ is of a repeller: by definition, the dynamical complement $S_{comp}$ admits an index pair of the form $(B, \emptyset)$, where $B$ is an isolating block for $S_{comp}$. Then for all $x\in B$, the forward trajectory through $x$ remains in $B$ because the exit set of $B$ is empty. Hence the $\omega$-limit $\omega(x)$ lies in the maximal invariant set $Inv(B)=S_{comp}$. \\
On the other hand, for all $x\in B\setminus{S_{comp}}$, there exists an entry-time $T\leq 0$ such that
$$ \forall t<T, \phi(t,x)\notin B.$$
Such trajectories accumulate on $S$, otherwise, they would entirely lie in $S_{comp}$ by definition of the dynamical complement. 

A similar argumentation proves the two  remaining claims.
\end{proof}

\begin{rem}
Of course, a very similar theorem holds if the set $S$ is isolated invariant:
\begin{enumerate}
\item If $h(S)$ is the index of a repeller, then there are trajectories through $x\notin S $ with $\alpha(x)\subset S$, but no trajectories through $x\notin S $ with $\omega(x)\cap S \neq \emptyset$.
\item If $h(S)$ is the index of an attractor, then there are trajectories through $x\notin S $ with $\omega(x)\subset S$, but no trajectories through $x\notin S $ with $\alpha(x)\cap S \neq \emptyset$.
\item If $h(S) $ is neither the index of an attractor nor the index of a repeller, then there are both trajectories  through $x\notin S $ with $\alpha(x)\subset S$ and through $y\notin S $ with $\omega(y)\subset S$
\end{enumerate}
\end{rem}

The following example illustrates how such theorems may be used.
\begin{ex} \label{trivial1}
We consider the  phase portrait as illustrated in Figure (\ref{ivialcompactempty}), which appears in the classification theorem of \cite{arllib}. This portrait shows a fixed point $S$ at infinity at the bottom of the phase portrait,  which is of isolated invariant complement $S_{comp}$.  The dynamical complement $S_{comp}$ in  Figure   (\ref{ivialcompactempty}) is of trivial Conley index, hence $\hat{h}(S)=\bar{0}$ too. This index is neither the index of a repeller nor the index of an attractor. Hence there are trajectories with $S$ in their $\alpha$-limit sets, and trajectories with $S$ in their $\omega$-limit sets, according to Theorem \ref{existence}. \\
Note that in the context of Conley index at infinity, a trivial index is no obstacle to the detection of connecting orbits. The connection structure will be analyzed  in Example (\ref{trivial2}).  
\begin{figure}
\begin{center}
\includegraphics[width=0.4\textwidth]{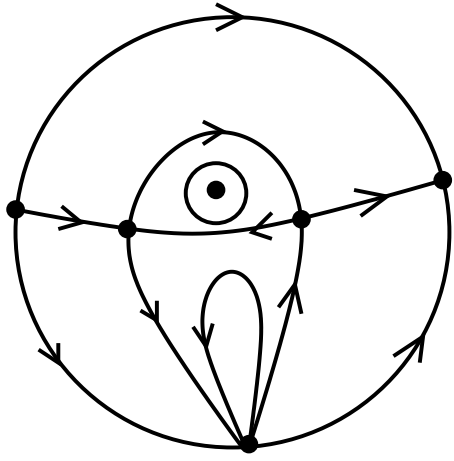}
\end{center}
\caption{ \label{ivialcompactempty}Phase portrait of Example (\ref{trivial1}) .} 
\end{figure}
\end{ex}

\section{Detection of Heteroclinics to Infinity} 
\subsection{Fundamental property}
The concept of isolated dynamical complement shows an interesting property on which our method for the detection of connection to some nonisolated invariant sets at infinity is based. Roughly speaking, this fundamental property says that leaving an isolating neighbourhood of $S_{comp}$ forces leaving all of them - leading to an accumulation on $S$. More precisely,  the following holds. 
\begin{prop} \label{exittoS}
Let $S\subset \mathcal{E}$ be an invariant set of isolated complement $S_{comp}$, and $N$ any isolating neighbourhood of the dynamical complement $S_{comp}$. Consider a trajectory $o(x)$ under the  flow $\phi$ on the Poincar\'e hemisphere $\mathcal{H}$ which connects to $S_{comp}$ in backward or forward time direction but leaves $N$ at a finite time $t_0$. Then this trajectory accumulates on $S$ in forward or backward time direction respectively. In other words, if $x\in N\setminus{S_{comp}}$ admits a $t_0$ such that $\phi(t_0, x)\notin N$, then:
$$
\begin{cases} 
\alpha(x)\subset S_{comp}\\
 \Rightarrow \omega(x)\cap S \neq \emptyset
\end{cases}
\begin{cases} 
\omega(x)\subset S_{comp}  \\
\Rightarrow\alpha(x)\cap S \neq \emptyset
\end{cases}
$$   
\end{prop}
\begin{proof}
Consider a growing sequence of isolating neighbourhoods $(N_n)_{n\in \mathbb{N}}$ of $S_{comp}$ with $\bigcap_{n\in \mathbb{N}}\mathcal{H}\setminus{N_n}=S$. Since $Inv(N)=S_{comp}=Inv(N_n)$ for every $n\in \mathbb{N}$, it holds for any $x\in \mathcal{H}$ that
$$
\phi(t_0, x) \notin N \Rightarrow \forall n \in \mathbb{N}, \exists t_n / \phi(t_n, x)\notin N_n.
$$
Since the sets $\mathcal{H}\setminus{N_n}$ shrink on $S$ as  $n$ grows,  $d(\phi(t_n,x), S)\rightarrow 0 $ as $n\rightarrow \infty$. \\
If in addition $\alpha(x)\subset S_{comp}$, then the sequence $(t_n)_{n \in \mathbb{N}}$ tends to $+\infty$, and, up to a subsequence, $\phi(t_n, x)$ converges to a point of $S$, so that $\omega(x)\cap S\neq \emptyset$.\\
Similarly, $\omega(x)\subset S_{comp}$ implies $\alpha(x)\cap S\neq \emptyset$. 
\end{proof}
\begin{rem}
Note that we can only guarantee a nonempty intersection of the $\omega$- or $\alpha$-limits with the invariant set $S$ of isolated invariant complement. An inclusion would not hold in general, since the $\omega$- or $\alpha$-limit could be a heteroclinic cycle containing equilibria in $S$ and equilibria in $S_{comp}$, for example, or even a homoclinic loop to $S$. 
\end{rem}

\subsection{Construction of an Ersatz Infinity} \label{constrersatz}

\begin{figure}
\begin{center}
\includegraphics[width=0.8\textwidth]{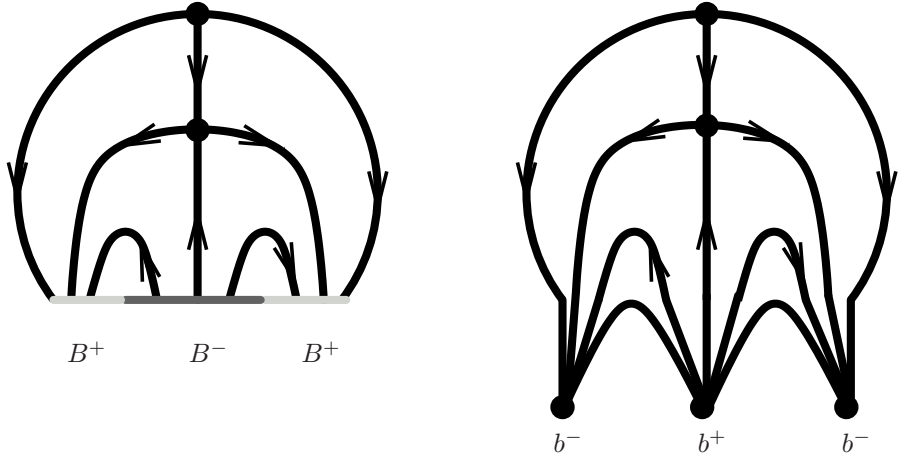}
\end{center}
\caption{Extended phase portrait for portrait 18 of Figure (\ref{classhamilton}).} 
\label{port18}
\end{figure}

We wish to prove the existence of connections between a given invariant set $S\subset \mathcal{E}$ at infinity with isolated invariant complement $S_{comp}$ and an isolated invariant set $Q\subset S_{comp}$. To this end, we replace $S$ with an isolated invariant ersatz for which Conley index techniques apply and use the fundamental property \ref{exittoS} stated above to conclude the existence of a connection to $S$. \\
We present here a construction of an ersatz which is realizable in many concrete cases. However, we do not pretend that it is optimal. It may be useful in some cases to choose  a more appropriate  ersatz infinity. The fundamental property allows us to be very flexible at this point: as long as we keep the flow unchanged in an isolating block for $S_{comp}$, we may complete the flow as we please outside of this block.\\
Let us give first an overview of the construction before we go into the technical details:
\begin{enumerate}
\item Find an isolating block $B$ for $S_{comp}$. Its boundary splits into immediate entrance and exit sets $\partial B = B^+  \cup B^-$, where $B^+ \cap B^-$ contains only a finite number of points of exterior tangency to $B$.
\item Choose retractions of $B^+$  and $B^-$ to sets $b^+$ and $b^-$ that are, typically,  less complicated (for example a union of single points).
\item Glue pieces of trajectories outside of $B$ that follow the retractions. The sets $b^{\pm} $ are isolated invariant, and repellers/attractors respectively.
\end{enumerate}
{\bf Construction of the extended phase space.}\\
Let us fix an invariant set $S$ in the sphere at infinity with isolated invariant complement $S_{comp}$. Furthermore, let $B$ be an isolating block for $S_{comp}$. There exist strong deformation retracts  $r_{\pm}: [0,1]\times B^{\pm}\rightarrow  B^{\pm}$ of the entrance and exit sets to subsets $ b^{\pm} \subset B^{\pm}$, i.e.
\begin{enumerate}
\item For all $x\in B^{\pm}$, $r^{\pm}(0,x)=x$.
\item  For all $x\in b^{\pm}$ and $s \in [0,1]$, $r^{\pm}(s,x)=x$.
\item For all $x\in B^{\pm}$, $r^{\pm}(1,x)\in b^{\pm}$
\end{enumerate}
Furthermore we require that the retractions $r^{\pm}$ be injective:  $$ \forall s\in[0,1[\text{, }\forall x_1, x_2 \in B^{\pm} \text{ , } r^{\pm}(s,x_1)=r^{\pm}(s,x_2) \Rightarrow x_1=x_2. \text{ (Inj) }$$ This injectivity property will be needed to construct a flow from the retracts without trajectories crossing each other. However this is no real restriction, as we see in the next lemma. 
\begin{lem}
Let $B$ be an isolating block for the isolated invariant dynamical complement $S_{comp}$ of $S\subset \mathcal{E}$. Then there exist strong deformation retracts satisfying the injectivity condition (Inj).
\end{lem}
\begin{proof}
Take $b^{\pm}=B^{\pm}$ and $r^{\pm}(s,.)=Id$.
\end{proof}
After we have fixed one pair of deformation retracts with property (Inj), we extend the block $B$ with subsets of $[0,1]\times B^{\pm}$ defined by the following:
\begin{equation} \label{betapm}
\beta^{\pm} :=   \{   (s,r^{\pm}(s,x)), s\in [0,1], x\in B^{\pm}     \}. 
\end{equation}
We now ``glue" the pieces $\beta^{\pm}$ on the boundary of $B$ by means of the equivalence relation 
$$\forall x,y\in B\cup \beta^+ \cup \beta^- , \text{ } x\sim y \text{  if }
\begin{cases}
 x=y \text{ or} \\
 x=(0,\xi)\in \beta^+\cup \beta^- \text{ and } y=\xi \in \partial B \text{ or}\\
 y=(0,\xi)\in \beta^+\cup \beta^- \text{ and } x=\xi \in \partial B.
 \end{cases}
 $$ 
 The extended phase space is defined as 
 \begin{equation} \label{hatH}
 \hat{\mathcal{H}}:= \frac{B\cup \beta^+ \cup \beta^-}{\sim} .
\end{equation}
For an illustration of the extension procedure on phase portrait 18 of Figure (\ref{classhamilton}), see Figure (\ref{port18}).
\vspace{0.4cm}\\
{\bf Construction of the extended flow} \\
We define a flow $\hat{\phi}$ on the extended phase space $\hat{\mathcal{H}}$ which follows the  flow $\phi$ on the block $B$, and follows the deformation retracts after leaving the block $B$. We first have to fix some notation. For each point $x\in B$, there exist entrance and exit times $T_-(x)\in [-\infty, 0]$,  $T_+(x)\in [0,+\infty]$, defined as follows:
\begin{eqnarray}
T_+(x):=\sup \{t\geqslant 0 / \phi([0,t],x)\subset B \}  \label{entextimes}, \\
T_-(x):=\inf \{t\leqslant 0 / \phi([t,0],x)\subset B \}. 
\end{eqnarray} 
Note that $T_{\pm}(x)$ may be infinite. In case they are finite,   $$\phi(T_-(x),x)\in B^+,\text{ and}$$  $$\phi(T_+(x),x)\in B^-.$$ The maps $T_{\pm}: B \rightarrow \mathbb{R}\cup \infty$ are continuous because $B$ is an isolating block. \\
On the other hand, if we consider a  point $(s,y)\in \beta^{\pm}$ where $s\in [0,1[$, then there is a point $x\in B^{\pm}$ such that $y=r^{\pm}(s,x)$, and this point is unique by virtue of  the injectivity condition (Inj). To $(s,y)\in \beta^{\pm}$, $s\neq 1$, we associate the  unique $Y^{\pm}(y):=\left( r^{\pm}(s, .)  \right)^{-1}(y) $. \\

\begin{defi} \label{extflow}
For an initial condition $x\in B$, the extended flow $\hat{\phi}(t,x)$ is defined as follows:
\begin{equation} \label{hatphiforxinB}
\hat{\phi}(t,x) = \begin{cases}
\phi(t,x) \in B \text{ if  } t\in [T_-(x), T_+(x)] ,\\
\left(   1-e^{T_+(x)-t} , r^-(1-e^{T_+(x)-t}, \phi(T_+(x),x))   \right) \in \beta^-  \text{ for } t\geqslant T_+(x)\\
 \left( 1-e^{t-T_-(x)} ,  r^+(1-e^{t-T_-(x)}, \phi(T_-(x), x))   \right) \in \beta^+\text{ for } t\leqslant T_-(x)
\end{cases}
\end{equation}
In particular $T_-(x)=-\infty$ implies $\hat{\phi}(t,x)=\phi(t,x) \in B$ for all $t\leqslant 0$, and  $T_+(x)=+\infty$ implies $\hat{\phi}(t,x)=\phi(t,x)\in B$ for all $t\geqslant 0$. Therefore, $Inv_{\hat{\phi}}(B)=Inv_{\phi}(B)=S_{comp}$.

For an initial condition $(s,y) \in \beta^-  $, $s\in[0,1[$, the associated point $Y^-(y)\in B^-$ with $y=r^-(s, Y^-(y))$
 belongs to $B$ so that $\hat{\phi}(t, Y^-(y))$ is well defined through (\ref{hatphiforxinB}) for all $t\in \mathbb{R}$. Hence we can define $\hat{\phi}(t,(s,y))$ by the following:
\begin{equation} \label{hatphiforysinMi-}
\hat{\phi}(t,(s,y)) = \begin{cases}
(1,y) \text{ if  } s=1\\
\hat{\phi}(t-\log(1-s), Y^-(y)) \text{ if }s\in [0,1[  
\end{cases}
\end{equation}

For an initial condition $(s,y)\in \beta^+$, $s\in[0,1[$, the associated point $Y^+(y)\in B^+$ with $y=r^+(s, Y^+(y))$
 belongs to $B$ so that $\hat{\phi}(t, Y^+(y))$ is well defined through (\ref{hatphiforxinB}) for all $t\in \mathbb{R}$. Hence we can define $\hat{\phi}(t,(s,x))$ by the following:
\begin{equation} \label{hatphiforysinMi+}
\hat{\phi}(t,(s,y)) = \begin{cases}
(1,y) \text{ if  } s=1\\
\hat{\phi}(\log(1-s)-t, Y^+(y)) \text{ if }s\in [0,1[  
\end{cases}
\end{equation}
\end{defi}
Again Figure (\ref{port18}) shows how the flow is extended for the phase Portrait 18 of Figure (\ref{classhamilton}).\\
The construction of the map $\hat{\phi}$ guarantees the following:
\begin{prop}
The map $\hat{\phi}:\mathbb{R} \times \hat{\mathcal{H}}\rightarrow \hat{\mathcal{H}}$ defined by \ref{extflow} is a continuous flow. 
In other words, it fulfills the following properties:
\begin{enumerate}
\item The map $\hat{\phi}$ is continuous.
\item $\hat{\phi}(0,.)=id_{\hat{\mathcal{H}}}$
\item For all $t,\tau \in \mathbb{R}$ and $p\in \hat{\mathcal{H}}$, $\hat{\phi}(t,\hat{\phi}(\tau,p))=\hat{\phi}(t+\tau,p)$.
\end{enumerate}

\end{prop}
\begin{rem}
\label{aubout}
By abuse of notation, we denote by $b^{\pm}$ the sets $$b^{\pm}=\big\{ \left(1,r^{\pm}(1,x)\right), x\in B^{\pm} \big\} \subset \hat{\mathcal{H}},$$
and call them ersatz infinities, because they replace the invariant set $S\subset \mathcal{E}$ and are isolated invariant so that the Conley index is able to deal with them. By Definition (\ref{extflow}) of the extended flow, they consist only of equilibria. Hence 
$$ \forall x \in B^+\text{, } \alpha_{\hat{\phi}}(x)=\lim_{t\rightarrow - \infty} \hat{\phi}(t,x)=\left( 1, r^+(1,x) \right) \in b^+ ,$$
$$ \forall x \in B^-\text{,  } \omega_{\hat{\phi}}(x)=\lim_{t\rightarrow \infty} \hat{\phi}(t,x)=\left( 1, r^-(1,x) \right) \in b^- .$$
In particular, $b^+$ is a repeller, while $b^-$ is an attractor. 
\end{rem}

\begin{rem} {\bf Smoothing the extension}\\
After this construction, we end up with both a manifold and a flow which are continuous but not a priori smooth. The smoothness is desirable for technical parts of the Conley index theory: it provides the existence of blocks $B$ which are manifolds with boundaries, and in particular guarantees the equality between the homology of the quotient $H_*(\frac{B}{B^-})$ and the relative homology $H_*(B,B^-)$.  Therefore we should require the smoothing of both $ \hat{\mathcal{H}}$ and $\hat{\phi}$. 
We conjecture that requiring the existence of a deformation smoothing both the extended phase space and the extended flow is no real restriction and is generically fulfilled. This is at least the case in the examples we show here, although we do not stress  this technical point for clarity.
\end{rem}
\subsection{Exhibit heteroclinics to infinity} 
We have now prepared all the ingredients for the detection of heteroclinic orbits to an invariant set $S$ at infinity of isolated invariant complement $S_{comp}$. There are many ways to formulate the detection of heteroclinics, for example in terms of connection matrices which are a powerful tool used to analyze the heteroclinic structure of a whole Morse decomposition. For simplicity, let us  rather stick to the cornerstone of this theory here. For this, let us consider  attractor repeller pairs in the extended phase space. More precisely, we assume that the extended phase space $\hat{\mathcal{H}}$ contains an invariant set $\hat{I}$ (under the extended flow $\hat{\phi}$) decomposable into an attractor repeller pair  $A,R\subset \hat{I}$ with the following characteristics:
\begin{enumerate}
\item One of the sets $A$, $R$ is an invariant subset of $b^- \cup b^+$. 
\item The other one is an invariant subset of $S_{comp}$. 
\end{enumerate}
In this situation, the ``nonidentity" $h(\hat{I})\neq h(A)\vee h(R)$ proves the existence of a heteroclinic connection $R\rightarrow A$ in the extended phase portrait $(\hat{\mathcal{H}}, \hat{\phi})$. As a consequence, there is an orbit leaving/entering the isolating block $B$ of $S_{comp}$ -- under both the flows $\hat{\phi}$ and $\phi$ because they coincide in $B$. By virtue of  Proposition \ref{exittoS}, this orbit  accumulates on $S$ in the phase portrait $(\mathcal{H}, \phi)$. This proves the following theorem.

\begin{theo}
\label{detect}
Let $\phi$ be a flow on the Poincar\'e hemisphere $\mathcal{H}$.
Consider an invariant set $S\subset \mathcal{E}$ of isolated invariant complement $S_{comp}$. Fix an isolating block $B$ for $S_{comp}$ and proceed to the construction of the extended phase space $\hat{\mathcal{H}}$ defined by (\ref{hatH}) and the extended flow $\hat{\phi}$ on $\hat{\mathcal{H}}$ defined by (\ref{extflow}). \\
Assume that there exist isolated invariant sets  $R$   and $ \hat{I}$ under the extended flow $\hat{\phi}$ with 
$$R \subset S_{comp},$$
$$ \hat{I} \supset R\cup b^-,$$ 
such that $(b^-, R)$ builds an attractor repeller decomposition of $\hat{I}$, and suppose further that the Conley index detects a heteroclinic connection $R\rightarrow b^-$ by means of 
$$h(\hat{I})\neq h(R)\vee h(b^-).$$
Then there is a heteroclinic connection $R\rightarrow S$ in the original flow $\phi$ on $\mathcal{H}$, or more precisely, there exists $x \in B$ with $\begin{cases} \alpha(x)\subset R \\ \omega(x)\cap S \neq \emptyset \end{cases}$.  \\
Analogously, assume that there exists isolated invariant sets  $A$   and $ \hat{I}$ under the extended flow $\hat{\phi}$ with 
$$A \subset S_{comp},$$
$$ \hat{I} \supset A\cup b^+,$$ 
such that $(A, b^+)$ builds an attractor repeller decomposition of $\hat{I}$, and suppose further that the Conley index detects a heteroclinic connection $R\rightarrow b^-$ by means of
$$h(\hat{I})\neq h(A)\vee h(b^+).$$
Then there is a heteroclinic connection $S\rightarrow A$ in the original flow $\phi$ on $\mathcal{H}$, or more precisely, there exists $x \in B$ with $\begin{cases} \omega(x)\subset A \\ \alpha(x)\cap S \neq \emptyset \end{cases}$.  
\end{theo}

\begin{ex} \label{trivial2}
Let us proceed to the extension for the Example (\ref{trivial1}) whose phase portrait is illustrated in Figure (\ref{ivialcompactempty}). 
\begin{figure}
\begin{center}
\includegraphics[width=0.4\textwidth]{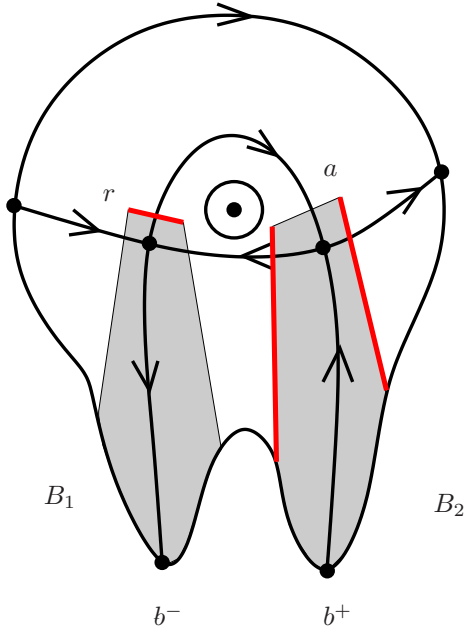}
\end{center}
\caption{ \label{blocksfortrivial} Extension and blocks for Example (\ref{trivial2}).} 
\end{figure}
The bottom equilibrium $S$ is invariant of isolated invariant dynamical complement $S_{comp}$. We choose an isolating block for $S_{comp}$ by cutting straight through the accumulation of homoclinic loops. The Figure (\ref{blocksfortrivial}) shows the extended phase space for this example together with isolating blocks $B_{1,2}$ for two attractor repeller pairs involving the ersatz infinity $b^+$ and $b^-$. \\
For the pair $(b^-,r)$, block $B_1$, the indices read:
$$
\begin{array}{lll}
\begin{cases}
h(\hat{\mathcal{H}}, \partial \hat{\mathcal{H}}; Inv(B_1))= \Sigma^1 \\
h(\hat{\mathcal{H}}; Inv(B_1))= \bar{0}\\
h( \partial \hat{\mathcal{H}}; Inv(B_1))= \Sigma^0
\end{cases} &
\begin{cases}
h(\hat{\mathcal{H}}, \partial \hat{\mathcal{H}}; r)= \Sigma^1 \\
h(\hat{\mathcal{H}}; r)= \Sigma^1\\
h( \partial \hat{\mathcal{H}}; r)= \bar{0}
\end{cases} &
\begin{cases}
h(\hat{\mathcal{H}}, \partial \hat{\mathcal{H}}; b^-)= \bar{0}\\
h(\hat{\mathcal{H}}; b^-)= \Sigma^0 \\
h( \partial \hat{\mathcal{H}}; b^-)= \Sigma^0
\end{cases} 
\end{array}
$$
As $ h(\hat{\mathcal{H}}; Inv(B_1))= \bar{0}\neq   \Sigma^1\vee \Sigma^0 =h(\hat{\mathcal{H}}; r) \vee h(\hat{\mathcal{H}}; b^-) $, we have by Theorem (\ref{detect}) a heteroclinic orbit $r\rightarrow S$
\vspace{0.3cm}\\
For the pair $(a,b^+)$, block $B_2$, the indices read:
$$
\begin{array}{lll}
\begin{cases}
h(\hat{\mathcal{H}}, \partial \hat{\mathcal{H}}; Inv(B_2))= \bar{0} \\
h(\hat{\mathcal{H}}; Inv(B_2))= \Sigma^1\\
h( \partial \hat{\mathcal{H}}; Inv(B_2))= \Sigma^1
\end{cases} &
\begin{cases}
h(\hat{\mathcal{H}}, \partial \hat{\mathcal{H}}; a)= \Sigma^1 \\
h(\hat{\mathcal{H}}; a)= \Sigma^1\\
h( \partial \hat{\mathcal{H}}; a)= \bar{0}
\end{cases} &
\begin{cases}
h(\hat{\mathcal{H}}, \partial \hat{\mathcal{H}}; b^+)= \Sigma^2\\
h(\hat{\mathcal{H}}; b^+)= \bar{0} \\
h( \partial \hat{\mathcal{H}}; b^+)= \Sigma^1
\end{cases} 
\end{array}
$$
As $ h(\hat{\mathcal{H}}, \partial \hat{\mathcal{H}}; Inv(B_2))= \bar{0} \neq   \Sigma^2\vee \Sigma^1 =h(\hat{\mathcal{H}}, \partial \hat{\mathcal{H}}; b^+) \vee h(\hat{\mathcal{H}} , \partial \hat{\mathcal{H}}; a) $, we have by Theorem (\ref{detect}) a heteroclinic orbit $S\rightarrow a$
\end{ex}

\begin{ex}
Let us exhibit heteroclinic connections to the bottom equilibrium in phase Portrait 18 of Figure (\ref{classhamilton}).
 The Figure (\ref{blocks18}) shows the extended phase space for this example together with isolating blocks $B_{1,2}$ for two attractor repeller pairs involving the ersatz infinity $b^+$ and $b^-$. \\
 \begin{figure}
\begin{center}
\includegraphics[width=0.8\textwidth]{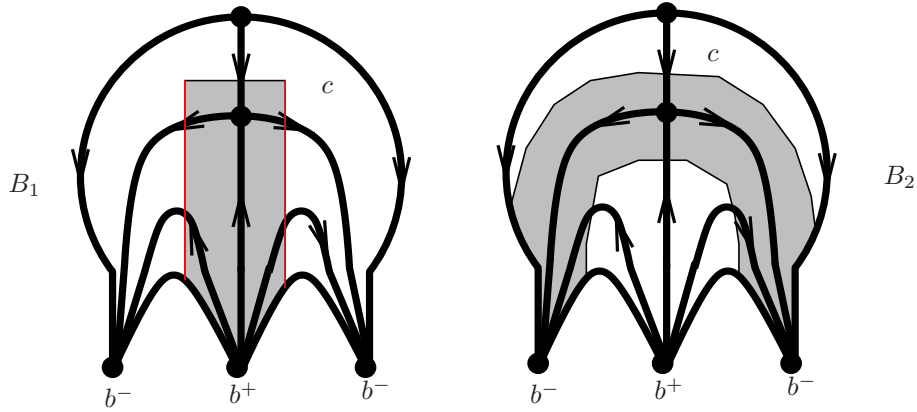}
\end{center}
\caption{ \label{blocks18} Extension and blocks for Portrait 18, Figure (\ref{classhamilton}).} 
\end{figure}
For the pair $(c,b^+)$, block $B_1$, the indices read:
$$
\begin{array}{lll}
\begin{cases}
h(\hat{\mathcal{H}}, \partial \hat{\mathcal{H}}; Inv(B_1))= \bar{0} \\
h(\hat{\mathcal{H}}; Inv(B_1))= \Sigma^1\\
h( \partial \hat{\mathcal{H}}; Inv(B_1))= \Sigma^1
\end{cases} &
\begin{cases}
h(\hat{\mathcal{H}}, \partial \hat{\mathcal{H}}; c)= \Sigma^1 \\
h(\hat{\mathcal{H}}; c)= \Sigma^1\\
h( \partial \hat{\mathcal{H}}; c)= \bar{0}
\end{cases} &
\begin{cases}
h(\hat{\mathcal{H}}, \partial \hat{\mathcal{H}}; b^+)= \Sigma^2\\
h(\hat{\mathcal{H}}; b^+)= \bar{0} \\
h( \partial \hat{\mathcal{H}}; b^+)= \Sigma^1
\end{cases} 
\end{array}
$$
As $ h(\hat{\mathcal{H}}, \partial \hat{\mathcal{H}}; Inv(B_1))= \bar{0} \neq   \Sigma^2\vee \Sigma^1 =h(\hat{\mathcal{H}}, \partial \hat{\mathcal{H}}; b^+) \vee h(\hat{\mathcal{H}} , \partial \hat{\mathcal{H}}; c) $, we have by Theorem (\ref{detect}) a heteroclinic orbit $S\rightarrow c$
\vspace{0.3cm}\\
For the pair $(b^-,c)$, block $B_2$, the indices read:
$$
\begin{array}{ll}
\begin{cases}
h(\hat{\mathcal{H}}, \partial \hat{\mathcal{H}}; Inv(B_2))= \Sigma^1 \\
h(\hat{\mathcal{H}}; Inv(B_2))= \Sigma^0\\
h( \partial \hat{\mathcal{H}}; Inv(B_2))= \Sigma^0 \vee \Sigma^0
\end{cases} &
\begin{cases}
h(\hat{\mathcal{H}}, \partial \hat{\mathcal{H}}; b^-)=\bar{0}\\
h(\hat{\mathcal{H}}; b^-)=  \Sigma^0 \vee \Sigma^0 \\
h( \partial \hat{\mathcal{H}}; b^-)=  \Sigma^0 \vee \Sigma^0 
\end{cases} 
\end{array}
$$
As $ h(\hat{\mathcal{H}}; Inv(B_2))= \Sigma^0 \neq   \Sigma^0\vee \Sigma^0 \vee \Sigma^1=h(\hat{\mathcal{H}}; c) \vee h(\hat{\mathcal{H}}; b^-) $, we have by Theorem (\ref{detect}) a heteroclinic orbit $c\rightarrow S$

\end{ex}

\section{Conclusions}
\subsection{Properties and Limits}
The construction of the extended phase portrait $(\hat{\mathcal{H}}, \hat{\phi})$ depends on the choice of the isolating block $B$ of the dynamical complement $S_{comp}$. In particular, the choice of $B$ influences the  topology of the ersatz infinity $b^{\pm}$, hence the ability of the Conley index to detect heteroclinic connections. Nevertheless it is sufficient to find one block $B$ for which Theorem (\ref{detect}) applies and detects a connection to the set $S$ of  isolated invariant dynamical complement $S_{comp}$.\\
The dependence of  the topology of $b^{\pm}$ on the choice of $B$ forbids us from defining the Conley index at infinity of $S$ as the sum $h_{\hat{\phi}}(b^+)\vee h_{\hat{\phi}}(b^-)$: this sounds like a good option but does not provide  a well defined index $\hat{h}(S)$.

The connections to the ersatz infinity $b^{\pm}$ occur in a one way street; this was not the case for the invariant $S$ itself, and this is the reason why it is not directly accessible to classical Conley index techniques. 
\vspace{0.5cm}\\
We described in the previous section how to detect heteroclinic orbits in the context of an attractor-repeller decomposition. More generally, Conley index and connection matrix theory allow one to analyze connection structure of a Morse decomposition involving more than two sets. Therefore it is  natural  to ask if we are able to deal with finitely many pairwise disjoint invariant sets $(S_i)_{1\leqslant i \leqslant n}\subset \mathcal{E}$ of isolated invariant complements $(S_{i_{comp}})_{1\leqslant i \leqslant n}$. The rough idea is to build a Morse decomposition in the extended phase space, analyze this decomposition via connection matrices and conclude  the existence of heteroclinic orbits in the original flow on the Poincar\'e hemisphere $\mathcal{H}$. We do not want to go into the details of the proofs here and refer the reader to \cite{myDiss}, Proposition (3.5.39)  and below. 
\begin{prop}
\label{collisoneigh}
Consider a flow $\phi$ on the Poincar\'e hemisphere $\mathcal{H}$, admitting a finite collection of pairwise disjoint invariant sets $(S_i)_{1\leqslant i \leqslant n}\subset \mathcal{E}$ of isolated invariant complements $(S_{i_{comp}})_{1\leqslant i \leqslant n}$.\\
Let $(K_{i_{\lambda}})_{1\leqslant i \leqslant n}$ be a collection of continuous families of isolating neighborhoods of the $(S_{i_{comp}})_{1\leqslant i \leqslant n}$. as in Proposition (\ref{def2}). Then 
\begin{enumerate}
\item The disjoint union $\bigcup_{i=1}^n \subset \mathcal{E}$ is an invariant set of isolated invariant dynamical complement $$\bigcap_{i=1}^n S_{i_{comp}}.$$
\item There exists $\lambda_0 \in [0,1[$ such that  $K_{\lambda}:=cl_{\mathcal{H}}\left( \bigcap_{i=1}^n int_{\mathcal{H}} (N_{i_{\lambda}}) \right)$, $\lambda\geqslant \lambda_0$,  forms  a continuous family of isolating neighborhoods for  $\bigcap_{i=1}^n S_{i_{comp}}$ in the sense of Proposition (\ref{def2}). 
\end{enumerate}
\end{prop}
To proceed to the construction of the extended phase portrait $(\hat{\mathcal{H}}, \hat{\phi})$, we need an isolating block $B$ for $\bigcap_{i=1}^n S_{i_{comp}}$. This block $B$ should fulfill some extra conditions allowing us to analyze the connection structure to the individual $S_i$'s via their ersatz $b_i^{\pm}$'s. More precisely we need that the immediate entrance/exit sets $B^{\pm}$ of $B$ split into $B_i^{\pm}$ leading after extension to ersatz $b_i^{\pm}$ of $S_i$.  For this purpose, we require isolating blocks $B_i$ of $S_{i_{comp}}$ satisfying the following condition:
\begin{equation}\label{nic}
\forall i\neq j \in \{ 1, \cdots , n \}, cl_{\mathcal{H}}(B_i^c) \cap cl_{\mathcal{H}}(B_j^c)=\emptyset 
\end{equation}
Condition (\ref{nic}) means that the blocks $B_i$ and $B_j$ are big enough so that their complements $B_i^c$ and $B_j^c$, neighborhoods of $S_i$ and $S_j$ respectively, are well separated. The existence of ``arbitrarily small" isolating blocks of $S_{i_{comp}}, S_{j_{comp}}$ is surely true, but not of``arbitrarily large" ones in general. Therefore the following  
\begin{conjecture}
\label{conj}
As soon as the $S_i$'s are pairwise disjoint invariant sets of isolated invariant complements, it is possible to find a collection of isolating blocks $(B_i)_{i\in \{ 1, \cdots , n \}}$ satisfying the condition (\ref{nic}).
\end{conjecture}
Once we have a collection $(B_i)_{i\in \{ 1, \cdots, n \}}$ of isolating blocks for $(S_{i_{comp}})_{i\in \{ 1, \cdots, n \}}$ satisfying condition (\ref{nic}), it follows that $$ B:=\bigcap_{i=1}^{n} B_i $$
is an isolating block for $ \bigcap_{i=1}^{n} S_{i_{comp}} $. We build the extended phase portrait $(\hat{\mathcal{H}}, \hat{\phi})$ on this block $B$. Suppose the extended flow $\hat{\phi}$ admits a Morse decomposition and classical connection matrix theory exhibits connections between an isolated invariant set $M\subset \bigcap_{i=1}^n S_{i_{comp}}$ and an ersatz infinity $b_{i_0}^{\pm}$. As a consequence, there exists a trajectory leaving/entering the block $B$ through its boundary $\partial_{\mathcal{H}}B=\bigcup_{i=1}^n \partial_{\mathcal{H}} B_i$, where the union is disjoint by virtue of  condition (\ref{nic}). Hence this trajectory leaves/enters also the block $B_{i_0}$ through its boundary. By Proposition (\ref{exittoS}), we know that the original trajectory under the flow $\phi$ on $\mathcal{H}$ accumulates on $S_{i_0}$. This shows that our method to detect heteroclinics is compatible  with connection matrix theory. Should the Conjecture (\ref{conj}) fail, the construction of the extended phase space should be performed for each $B_i$ individually: The truth of the Conjecture (\ref{conj}) does not affect the main Theorem (\ref{detect}), only the efficiency of its application.

As noticed in Remark (\ref{aubout}), the ersatz $b_i^+$ are repellers, the $b_i^-$ are attractors. They are hence maxima/minima  of  chains of the Morse ordering. 

Let us again consider the case of a single invariant set $S$ of isolated invariant complement in order to describe the type of information that leads to the detection of heteroclinic  connections. Note that since  $b^+$ is a repeller, its homological Conley index  $H_*(b^+)$ has a nontrivial $H_n(b^+)$, where $n$ is the dimension of our original phase space, say $n=\dim(\mathcal{H})$. A nontrivial connection map $\delta: H_n(b^+)\rightarrow H_{n-1}(A)$ for an isolated invariant $A\subset S_{comp}$ detects a  heteroclinic connection $b^+ \rightarrow A$ in $\hat{\mathcal{H}}$, hence $S\rightarrow A$ in $\mathcal{H}$ by Theorem (\ref{detect}). Furthermore, the topology of $b^+$ itself may produce other nontrivial homologies $H_j(b^+)$, so that other connection maps $\delta: H_j(b^+)\rightarrow H_{j-1}(A')$ may detect other heteroclinics  $S\rightarrow A'$. \\
Analogously, the ersatz $b^-$ is an attractor, its homological Conley index on the 0-th level is nontrivial so that a nontrivial connection map $\delta: H_1(R)\rightarrow H_0(b^-)$ detects a connection $R\rightarrow S$. Further nontrivial  homology groups $H_j(b^-)$ of the homological Conley index of $b^-$ and nontrivial connection maps $\delta: H_{j+1}(R')\rightarrow H_j(b^-)$ may detect further heteroclinic connections  $R' \rightarrow S$. 
 
\subsection{Genericity Questions} \label{genericity}
The following theorem has been proved in \cite{gonzvel}. It deals with genericity of hyperbolic equilibria at infinity in the class of polynomial vector fields. 
\begin{theo}
\label{strucstab}
Let $\mathcal{X}$ be the set of polynomial vector fields on $\mathbb{R}$ of degree smaller than $d\in \mathbb{N}$. Furthermore for every $f\in \mathcal{X}$, let us denote by $f_{\mathcal{E}}$ the vector field induced by $f$ on the sphere at infinity $\mathcal{E}$ after compactification and normalization. The set 
$$ Hyp=\{ f \in \mathcal{X} / f_{\mathcal{E}}\text{  admits only hyperbolic equilibria on } \mathcal{E}  \}$$
is open and dense in $\mathcal{X}$.
\end{theo}
In other words the equilibria at infinity are structurally stable (in $\mathcal{X}$) as soon as they are hyperbolic in the sphere at infinity. An accumulation of homoclinic loops, if it takes place in the finite part $\mathcal{H}\setminus{\mathcal{E}}$ of the Poincar\'e hemisphere, is structurally stable. Therefore the field of application of our methods is not restricted to pathological cases.  
\subsection{Perspectives for PDEs}

The concepts we have developped in this paper are promising for the study of blow up phenomena in partial differential equations. 

Compactification has already been  used  to analyze the dynamics of the global unbounded attractor for the asymptotically linear equation $u_t=Au+g(u)$, where $u$ belongs to an Hilbert space, $g$ is sublinear, and $A$ nice enough. See \cite{nitsan, myDiss} for details. In this example, grow up takes place and all equilibria at infinity admits classical Conley index in the sense of \cite{ryba}. 

For  PDE's whose leading terms are  homogenous of degree $d\geqslant1$ (for instance $u^d$ or $u^{d-1}\Delta u$, ...) we do not expect classical Conley index theory to be able to analyze the connection structure of the global unbounded attractor. 

In order to be able to apply our methods, some requirements are needed. First of all, we need to reduce the problem to a finite dimensional situation. Let us briefly explain why. 
The classical Conley index has been adapted to infinite dimensional situations and succesfully applied to partial differential equations. We do not want to develop this theory here: for the details, see \cite{ryba}.  Let us illustrate our dilemma on a dynamical system in a infinite dimensional phase space, compactified on an infinite dimensional Poincar\'e hemisphere $\mathcal{H}$, admitting an invariant set $S$ at infinity of isolated invariant dynamical complement $S_{comp}$. Suppose further  that $S_{comp}$  is a hyperbolic equilibrium with finite dimensional unstable manifold of dimension $u$. The classical Conley index in the sense of \cite{ryba} is well defined:  
\begin{enumerate}
\item If the direction pointing into the Poincar\'e hemisphere is stable, the Conley index reads
$$h(S_{comp})= \begin{cases} 
                              h(\mathcal{H}; S_{comp})= \Sigma^u,    \\
                              h(\mathcal{H}, \mathcal{E}; S_{comp})= \bar{0}\\
                              h(\mathcal{E}; S_{comp})= \Sigma^u.
                              \end{cases}$$
\item If the direction pointing into the Poincar\'e hemisphere is unstable, the Conley index reads
$$h(S_{comp})= \begin{cases} 
                              h(\mathcal{H}; S_{comp})= \bar{0},    \\
                              h(\mathcal{H}, \mathcal{E}; S_{comp})= \Sigma^u\\
                              h(\mathcal{E}; S_{comp})= \Sigma^{u-1}.
                              \end{cases}$$
\end{enumerate}
The Conley index at infinity $\hat{h}(S)$ was defined by time duality: in this case it would not be defined in the sense of \cite{ryba} because $S_{comp}$ has infinite dimensional unstable manifold  under reversed time direction. Or, to see it in a different way, our definition applied formally would produce an index at infinity equal to an infinite dimensional pointed sphere of codimension $u$: such a sphere is homotopy equivalent to a point, i.e.  a trivial index, which  does not make any sense. The same can be said about the index of the erstaz infinity $b^+$ which is a repeller: the Conley index of $b^+$ under the extended flow $\hat{\phi}$,  $h_{\hat{\phi}}(b^+)$, is not well defined, or absurdly trivial. Hence theorem \ref{detect} cannot be applied directly in infinite dimensional context.

For these reasons, we have to reduce the problem to a finite dimensional one. A good way of doing so is to prove the existence of an invariant manifold -- the inertial manifold -- containing the global unbounded attractor. See \cite{nitsan} for such an approach: in this work, grow up phenomena are studied, and the connection structure of the attractor is analyzed via a $y$-map type of tool.\\
An alternative approach is to use the Conley Index introduced by Izydorek and Rybakowski in \cite{IzyRyba}: the Conley Index at infinity $\hat{h}(S)$ of the example above would make sense in the finite dimensional Galerkin approximations, as well as the indices of the ersatz infinities. Hence exhibiting heteroclinic connections for  each finite dimensional Galerkin approximations via the extended flow would allow to conclude to the existence of heteroclinic connections in the original PDE.

\bibliography{litterature}{}

\begin{thebibliography}{10}

\bibitem{ALGM}
A.A. Andronov, E.A. Leontovich, I.I. Gordon, and A.G. Maier.
\newblock {\em Qualitative Theory of Second-order Dynamic Systems}.
\newblock {Halsted Press Book}, 1973.

\bibitem{arllib}
Joan~C. Art\'es and Jaume Llibre.
\newblock Quadratic hamiltonian vector fields.
\newblock {\em Journal of Differential Equations}, 1994.

\bibitem{nitsan}
Nitsan Ben-Gal.
\newblock {\em Grow-up Solutions and Heteroclinics to Infinity for Scalar
  Parabolic PDE's}.
\newblock PhD thesis, Brown University, Providence, USA, 2009.

\bibitem{conley}
Charles Conley.
\newblock Isolated invariant set and the {M}orse index.
\newblock {\em CMBS Regional Conf.~Ser.~Math.}, 38, 1978.

\bibitem{gonzvel}
Enrique~A. Gonz\'alez~Velasco.
\newblock Generic properties of polynomial vector fields at infinity.
\newblock {\em Trans. Am. Math. Soc.}, 143:201--222, 1969.

\bibitem{myDiss}
Juliette Hell.
\newblock {\em Conley Index at Infinity}.
\newblock PhD thesis, Freie Universit\"at, Berlin, Germany, 2010.

\bibitem{IzyRyba}
Marek Izydorek and Krzysztof~P. Rybakowski.
\newblock On the conley index in hilbert spaces in the absence of uniqueness.
\newblock {\em Fundamenta Mathematicae}, 171, 2002.

\bibitem{kappos3}
Efthimios Kappos.
\newblock The {C}onley index and global bifurcations. {I}: Concepts and theory.
\newblock {\em Int. J. Bifurcation Chaos Appl. Sci. Eng.}, 5(4):937--953, 1995.

\bibitem{kappos4}
Efthimios Kappos.
\newblock The {C}onley index and global bifurcations. {II}: Illustrative
  applications.
\newblock {\em Int. J. Bifurcation Chaos Appl. Sci. Eng.}, 6(12b):2491--2505,
  1996.

\bibitem{kappos2}
Efthimios Kappos.
\newblock Compactified dynamics and peaking.
\newblock In {\em Circuits and Systems, 2000. Proceedings. ISCAS 2000 Geneva.},
  2000.

\bibitem{mccord}
Christopher McCord.
\newblock Poincar\'e--{L}efschetz duality for the homology {C}onley index.
\newblock {\em Transactions of the American Mathematical Society}, 1992.

\bibitem{messias}
Marcelo Messias.
\newblock Dynamics at infinity and the existence of singularly degenerate
  heteroclinic cycles in the {L}orenz system.
\newblock {\em J.Phys. A: Math Theor}, 2009.

\bibitem{mrozek}
Marian Mrozek and Roman Srzednicki.
\newblock On time-duality of the {C}onley index.
\newblock {\em Result. Math.}, 24(1-2):161--167, 1993.

\bibitem{perko}
Lawrence Perko.
\newblock {\em Differential Equations and Dynamical Systems}.
\newblock Springer-Verlag, 1991.

\bibitem{ryba}
Krzystof~P. Rybakowski.
\newblock {\em The Homotopy Index and Partial Differential Equations}.
\newblock Springer-Verlag., 1987.

\bibitem{schlovu}
Dana Schlomiuk and Nicolae Vulpe.
\newblock Geometry of quadratic differential systems in the neighborhood of
  infinity.
\newblock {\em J. Differ. Equations}, 215(2):357--400, 2005.

\bibitem{smollerbook}
Joel Smoller.
\newblock {\em Shock Waves and Reaction-diffusion Equations.}
\newblock { Springer-Verlag}, 1983.

\end{thebibliography}
\bibliographystyle{plain}
\end{document}